\documentclass[final,leqno]{siamltex}
\usepackage{amstext}
\usepackage{graphicx,subfigure}
\usepackage{amsfonts}
\usepackage{amsmath}
\usepackage[notref,notcite]{showkeys}

\newcommand\cO{{\cal O}}

\renewcommand\epsilon{\varepsilon}
\newcommand\eps{\varepsilon}
\renewcommand\phi{\varphi}
\newcommand{\pp}{\partial}

\newcommand{\tr}{\tilde{r}}
\newcommand{\ttt}{\tilde{t}}
\newcommand{\sr}[1][]{\sigma_{#1r}}
\newcommand{\si}[1][]{\sigma_{#1i}}
\newcommand{\lr}{\lambda_r}
\newcommand{\li}{\lambda_i}
\newcommand{\dsr}[1][]{\dot{\sigma}_{#1r}}
\newcommand{\dsi}[1][]{\dot{\sigma}_{#1i}}
\newcommand{\R }{\bf R}
\newcommand{\re}{\mathbb{R}}
\newcommand{\tC}{\widetilde{C}}

\begin{document}

\baselineskip=12pt
\bibliographystyle{siam}

\title{(In-)Stability of singular equivariant solutions to the Landau-Lifshitz-Gilbert equation}

\author{Jan Bouwe van den Berg\thanks{Dept.\ of Mathematics, VU University Amsterdam, de Boelelaan 1081, 1081 HV Amsterdam, the Netherlands, {\tt janbouwe@math.vu.nl}} \and J.F. Williams\thanks{Dept.\ of Mathematics, Simon Fraser University, 
Burnaby, Canada, {\tt jfw@math.sfu.ca}}}

\date{\today}

\maketitle

\begin{abstract}
In this paper we use formal asymptotic arguments to understand the
stability properties of equivariant solutions to the
Landau-Lifshitz-Gilbert model for ferromagnets.  We also analyze both the
harmonic map heatflow and Schr\"odinger map flow limit cases.  All
asymptotic results are verified by detailed numerical experiments, as well as a robust topological argument.
The key result of this paper is that blowup solutions to these problems
are co-dimension one and hence both unstable and non-generic. 

Finite time blowup solutions are thus far only known to arise
in the harmonic map heatflow in the special case of radial symmetry.
Solutions permitted to deviate from this symmetry remain global for all time
but may, for suitable initial data, approach arbitrarily close to blowup.
A careful asymptotic analysis of solutions near blowup shows that finite-time blowup corresponds to a saddle fixed point in a low dimensional dynamical system.  Radial symmetry precludes motion anywhere but on the stable manifold towards blowup.

The Landau-Lifshitz-Gilbert problem is not invariant under radial symmetry. Nevertheless, a similar scenario emerges in the equivariant setting: blowup is unstable. To be more precise, blowup is co-dimension one both within the equivariant symmetry class and in the unrestricted class of initial data.
The value of the parameter in the Landau-Lifshitz-Gilbert equation plays a very subdued role in the analysis of equivariant blowup, leading to identical blowup rates and spatial scales for all parameter values. 
One notable exception is the angle between solution in inner scale (which bubbles off) and outer scale (which remains), which does depend on parameter values.

Analyzing near-blowup solutions, we find that in the inner scale these solution quickly rotate over an angle $\pi$. As a consequence, for the blowup solution it is natural to consider a continuation scenario after blowup where one immediately re-attaches a sphere (thus restoring the energy lost in blowup), yet rotated over an angle $\pi$. This continuation is natural since  it leads to continuous dependence on initial data.

\end{abstract}


\section{Introduction}
\setcounter{section}{1}
\setcounter{equation}{0}
In this paper we are interested in the existence and stability of
finite-time singularities of the Landau-Lifshitz-Gilbert equation for maps
from the unit disk (in the plane) to the surface of the unit sphere, $m: D^2 \to S^2$:
\begin{equation}
\label{LLG1}
\left\{
\begin{aligned}
\frac{\pp m}{\pp t} & = \alpha \, m \times \Delta m 
                            - \beta \, m\times(m\times \Delta m) , \\
m(x,t) &= m_b(x) \qquad |x| = 1, \\
m(x,0) & = m_0(x) .
\end{aligned}
\right .
\end{equation}
We will always require the damping term $\beta \ge 0$ and take $\alpha^2 +
\beta^2 = 1$ without loss of generality.  This problem preserves the length
of the vector $m$, i.e., $|m_0(x)| = 1$ for all~$x$ implies that 
$|m(x,t)| = 1$ for all positive time (for all $x$).

In the case $ \alpha \ne 0, \beta > 0$ this equation arises as a model for
the exchange interaction between magnetic moments in a magnetic spin system
on a square lattice \cite{HS, KP}.  Taking $\alpha = 0$ recovers the harmonic
map heatflow which is a model in nematic liquid crystal flow \cite{BBCH}.
It is also of much fundamental interest in differential geometry
\cite{Struwe}.  Finally, the conservative case $\beta = 0$ is the
Schr\"odinger map from the disk to the sphere, which is a model of current
study in geometry \cite{CSU,GustKangTsai,GKT2}.

Stationary solutions in all these cases are \emph{harmonic maps}. This allows us to analyze singularity formation in a
unified manner for all parameter values.  Traditionally, the Landau-Lifshitz-Gilbert equation is posed with
Neumann boundary conditions, but this does not affect the local structure of
singularities, should they arise.
We note that in the harmonic map literature the second term on the right-hand side of the differential equation~(\ref{LLG1}) is often rewritten using the identities
\[ m\times(m\times \Delta m)  = - \Delta m + (\Delta m, m) m =  - \Delta m - 
\| \nabla m \|^2 m, \]
where $(\cdot,\cdot)$ denotes the inner product in $\re^3$ and 
$\| \nabla m \|^2 = \sum_{i=1}^2\sum_{j=1}^3 (\partial_i w_j)^2$, with $w_i$ the components of $w$.

As discussed in much greater detail in Sections~\ref{sec:formulation} and ~\ref{sec:previousresults}, there are initial data for which the solution to~(\ref{LLG1}) becomes singular in finite time. In this paper we analyze this \emph{blowup} behaviour, in particular its \emph{stability} properties under (small) perturbations of the initial data. Considering initial data that lead to blowup, the question is whether or not solutions starting from slightly different initial data also blowup.
Our main conclusion is that blowup is an \emph{unstable co-dimension one} scenario. With this in mind we also investigate the behavior of solutions in ``near-miss'' of blowup, and the consequences this has for the continuation of the blowup solution after its blowup time.

\subsection{Problem formulation}
\label{sec:formulation}

We will consider two formulations for equation (\ref{LLG1}).  The first is
the so-called \emph{equivariant} case: using polar coordinates $(r,\psi)$ on the unit disk $D=D^2$, these are solutions of the form
\begin{equation}
\label{LLG_decomp1}
m(t,r,\psi) = \left(
\begin{array}{c}
\cos(n\psi) u(r,t) - \sin(n\psi) v(r,t)\\
\sin(n\psi) u(r,t) + \cos(n\psi) v(r,t) \\
w(r,t)
\end{array}
\right) , 
\end{equation}
which have the (intertwining) symmetry property
$ m(t,\cdot) \circ R_2^\omega = R_3^{n\omega} \circ m(t,\cdot)$ 
for all~$\omega$ and each fixed~$t$,
where $R_2^\omega$ is a rotation over angle $\omega$ around the origin in the plane~$\re^2$, while $R_3^\omega$ is a rotation over angle $\omega$ around the $z$-axis in $\re^3$.

The components $(u,v,w)$ then satisfy the pointwise constraint $u^2+v^2+w^2=1$, as well the differential equations
\begin{equation}
\label{LLG_3comp}
\left\{
\begin{aligned}
u_t &=  \alpha \left(v \Delta w - \left(\Delta v - \frac{n^2}{r^2} v \right) w \right) + 
\beta \left(\Delta u - \frac{n^2}{r^2}u + A u \right) , \\
v_t &= \alpha \left(-u \Delta w + \left(\Delta u - \frac{n^2}{r^2} u \right) w \right) + 
\beta \left(\Delta v - \frac{n^2}{r^2}v + A v \right) , \\
w_t & = \alpha \left(u \Delta v - v \Delta u \right)  + \beta \left(\Delta w + A w \right)  ,\\
\end{aligned}
\right.
\end{equation}
where
\[ \Delta = \frac{\pp^2}{\pp r^2} + \frac{1}{r}\frac{\pp}{\pp r} \qquad \mbox{and} \qquad
A \equiv u_r^2 + v_r^2 + w_r^2 + \frac{n^2}{r^2}\left(u^2 + v^2\right). \]
We will take $n = 1$ in what follows, except in Section~\ref{sec:AsympN}.

Alternatively, we can parametrize the solutions on the sphere via the Euler angles:
\begin{equation}
m(t,r,\psi) = \left( 
\begin{array}{c}
\cos [\psi+\phi(r,t)] \sin \theta(r,t) \\
\sin [\psi+\phi(r,t)] \sin \theta(r,t) \\
\cos \theta(r,t) 
\end{array}
\right ) ,
\label{LLG_2form}
\end{equation}
where the equations for $\theta$ and $\phi$ are given by 
\begin{equation}
\left\{
\begin{aligned}
\beta \theta_t + \alpha \sin \theta \phi_t  & = \theta_{rr} + \frac{1}{r} \theta_r - \frac{\sin 2\theta}{2}
\left(\frac{1}{r^2} + \phi_r^2\right) , \\
 \beta \phi_t - \frac{\alpha}{ \sin \theta} \theta_t  &= \phi_{rr} + \frac{1}{r} \phi_r + \frac{\sin 2\theta}{\sin^2 \theta} \phi_r \theta_r .
\end{aligned}
\right.
\label{LLG_2comp}
\end{equation}
We note that due to the splitting $\psi+\phi(r,t)$ in (\ref{LLG_2form})
the system (\ref{LLG_2comp}) has \emph{one} spatial variable.
In this equivariant case the image of one radius in the disk thus fixes the entire map (through rotation) and we write $m(t,r)=m(t,r,0)$.

In the special case $\alpha = 0$ and $\beta = 1$ only, there are \emph{radially symmetric} solutions of the form $\phi \equiv$~constant, reducing the system to a single equation
\begin{equation}
\label{radHM}
\theta_t = \theta_{rr} + \frac{1}{r}\theta_r - \frac{\sin2\theta}{2r^2} .
\end{equation}

\subsection{Previous results}
\label{sec:previousresults}

It is well known that not all strong solutions to the radially symmetric harmonic map heatflow (\ref{radHM}) are global in time. Equation (\ref{radHM}) is $\pi$-periodic in $\theta$. Supplemented with the (finite energy) boundary condition $\theta(0)=0$, it only has stationary solutions of the form
$u_q^{\infty} = 2 \arctan{q r}$
for $q \in \R$. Hence with prescribed boundary data $\theta(0,t) = 0$ and $\theta(1,t) = \theta^* > \pi$ there is no accessible stationary profile.  However, 
there is an associated Lyapunov functional,
\[
  E(t) =  \pi \int_0^1 \left[ \theta_r(t,r)^2  
                         + \frac{\sin^2 \theta(t,r)}{r^2} \right] r \, dr  
\]
whose only stationary points are the family $u_q^\infty$.  It is this paradox that leads to blowup: there is a finite collection of (possibly finite)  times at which $u(0,t)$ ``jumps'' from $n\pi$ to $(n\pm 1)\pi$, losing $4\pi$ of energy~\cite{Struwe,BertschHulshofHout}. The structure of the local solution close to the jumps (in time and space) is known, which allows us to analyze the stability of these solutions.

The fundamental result in this area is due to Struwe \cite{Struwe} who first showed that solutions of the harmonic map heatflow could exhibit the type of  jumps described above and derived what the local structure of the blowup profile is.  Chen, Ding and Ye \cite{CDY} then used super- and sub-solution arguments, applicable only to (\ref{radHM}), to show that finite-time blowup must occur when $u(0,t) = 0$ and $u(1,t) \ge \pi$.  The blowup rate and additional structural details were determined through a careful matched asymptotic analysis in \cite{BHK}.

The analysis is
based on the original result of Struwe who showed that any solution which
blows up in finite time must look locally (near the blowup point) like a rescaled harmonic map at the so-called quasi-stationary scale.
That is, there is a scale $r=O(R(t))$ on which the solution takes the form
\begin{equation}
\label{InnerProf}
\theta(t,r) \to 2 \arctan \left(\frac{r}{R(t)}\right) \qquad \mbox{where} \;\; R(t) \to 0 \quad \mbox{as} \; t \to T.
\end{equation}
From \cite{BHK} it is known that for generic initial data one has
\begin{equation}
\label{RadRate}
R(t) \sim \kappa \frac{T-t}{|\ln(T-t)|^2} \qquad \mbox{as}\quad t \to T
\end{equation}
for some $\kappa > 0$ and blowup time $T > 0$, which both depend on the initial data. This result is intriguing
as the blowup rate is very far from the similarity rate of $\sqrt{T-t}$~\cite{AngenentHulshofMatano}. 

While the blowup rate~(\ref{RadRate}) was derived in~\cite{BHK} for the harmonic map heatflow, i.e.\ $\alpha=0$, $\beta=1$, in this paper we demonstrate that formal asymptotics imply that this rate is universal for all parameter values of $\alpha$ and $\beta$. 
Very recently, proofs of the blowup rate~(\ref{RadRate}) have appeared for the harmonic map heatflow~\cite{RS} as well as the Schr\"odinger map flow~\cite{MRR}, i.e., the two limit cases of the Landau-Lifshitz-Gilbert problem.

It is common to consider radial symmetry when analyzing the blowup
dynamics of many reaction-diffusion equations.  Typically there one can
show that there must be blowup using radially symmetric arguments.
Moreover, numerical experiments generically show that rescaled solutions
approach radial symmetry as the blowup time is approached.

For the harmonic map problem the proof of blowup solutions due to Cheng,
Ding and Ye is completely dependent on the radial symmetry.  Moreover,
there are stationary solutions to the problem which are in the homotopy
class of the initial data, but which are not reachable under the radial symmetry constraint.  This begs the question:
\emph{What happens when we relax the constraint of radially symmetric
initial data?}

The above description is mainly restricted to the harmonic map problem ($\alpha=0$) which has received considerably more attention than the general Landau-Lifshitz-Gilbert equation.  
Before addressing the question of stability under non-radially symmetric constraints for the harmonic map heatflow problem, we first show that blowup solutions are still expected for the full Landau-Lifshitz-Gilbert equation with $\alpha, \beta > 0$, see Section~\ref{sec:Top}. 
We note that 
\[
  E(t) = \pi \int_0^1 \left[ \theta_r(t,r)^2  
    + \sin^2 \theta(t,r) \left( \phi_r(t,r)^2 +\frac{1}{r^2} \right)
            \right] r \, dr  
\]
is a Lyapunov functional for the equivariant problem (\ref{LLG_2comp}) as long as $\beta>0$, whereas it is a conserved quantity for the Schr\"odinger map flow ($\beta=0$).  

We discuss the question of stability for the full problem in a
uniform manner.
The topological argument in Section~\ref{sec:Top} suggests that blowup is co-dimension one, and this is indeed supported by the asymptotic analysis in Section~\ref{sec:Asymp} and the numerics in Section~\ref{sec:Num}.
The main quantitative and qualitative properties turn out to be independent of the parameter values, except for the angle between sphere that bubbles off and the remaining part of the solution.
In Section~\ref{sec:MoreAsymp} we 
analyze near-blowup solutions. These solution rotate quickly over an angle $\pi$ in the inner scale. For the blowup solution this implies a natural  continuation scenario (leading to continuous dependence on initial data) after the time of blowup: the lost energy is restored immediately by re-attaching a sphere, rotated over an angle $\pi$.
Finally, in Section~\ref{sec:AsympN} we present the generalization to the case $n \geq 2$, followed by a succinct conclusion in Section~\ref{sec:Concl}.

\section{The global topological picture}
\label{sec:Top}

We present a topological argument to corroborate that blowup is
co-dimension one. It does not distinguish between finite and infinite time
blowup. Since the argument relies on dissipation,
it works for $\beta>0$, but since the algebra is
essentially uniform in $\alpha$ and $\beta$, as we shall see in Section~\ref{sec:Asymp}, we would argue that the
situation for the Schr\"odinger map flow is the same.

Let us first consider the equivariant case, where, as explained in Section~\ref{sec:formulation}, the image of one radius in the disk fixes the entire map, and we write $m(t,r)=m(t,r,0)$.

Let $m(t,0)=N$ (the north pole) and $m(t,1) = m_b$.
In the notation using Euler angles from Section~\ref{sec:formulation}, 
by rotational symmetry we may assume that $m_b=(\theta_b,0)$, $\theta_b \in
[0,\pi]$. 
The \emph{only} equilibrium  configuration satisfying these boundary
conditions is $m = (\theta,\phi) = (2 \arctan q r ,0)$, where $q=\tan
(\theta_b/2)$. Note that for $\theta_b=\pi$ there is no equilibrium, hence
blowup must occur for all initial data in that case~\cite{AngenentHulshof}.

For $\theta_b \in [0,\pi)$, i.e.\ $m_b \neq S$, we shall construct a \emph{one 
parameter family} of initial data $m_0(r;s)$, and we argue that at least one of the corresponding solutions blows up. Since the presented argument is
topological, it is robust under perturbations, hence it proves that blowup is (at most) co-dimension one. The matched asymptotic analysis in Section~\ref{sec:Asymp} confirms this co-dimension
one nature of blowup.

\begin{figure}[t]
\centerline{\includegraphics[height=4cm]{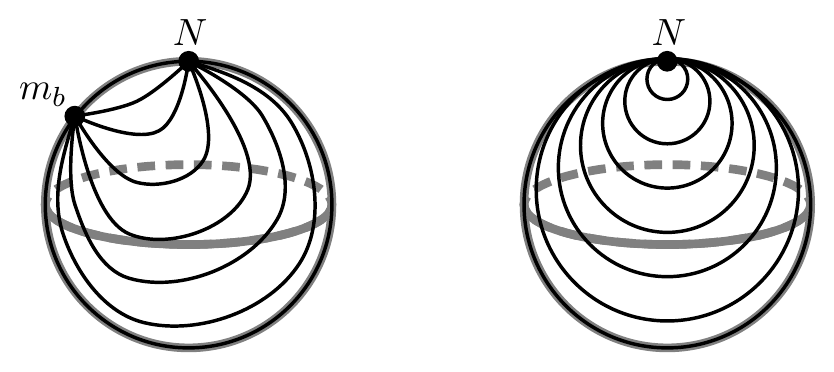}}
\caption{One parameter families of initial data for the equivariant case;
several members of half of each family are shown (the other half lives on the hemisphere facing away from us).
Left: for $m_b \neq N$ one may obtain such a family for example by
stereographic projection (w.r.t.~$N$) of all straight lines through the
point in the plane corresponding to $m_b$. Right: for $m_b = N$ one can 
choose the stereographic projection of parallel lines covering the plane.}
\label{f:cover1}
\end{figure}

We choose one-parameter families of initial data as follows.
The family of initial data will be parametrized by $s \in S^1$, or $[0,1]$ with the end points identified. Let $m_0(r;s)$ be a continuous map from
$[0,1] \times S^1$ to $S^2$, such that
\begin{equation}
	\label{e:endpoints}
   m_0(0;s)=N  \quad\text{and}\quad  
   m_0(1;s)=m_b \qquad\text{for all }s.
\end{equation}
We may then view $m_0$ as a map from $S^2 \to S^2$ by identifying $\{0\} \times S^1$ and $\{1\} \times S^1$ to points.
Now choose any continuous family $m_0(r;s)$ satisfying (\ref{e:endpoints}) such that it represent a \emph{degree 1} map from $S^2$ to itself.
One such choice is obtained by using the stereographic projection
\[
  T(x,y) = \left(
   \frac{2x}{1+x^2+y^2},\frac{2y}{1+x^2+y^2},\frac{-1+x^2+y^2}{1+x^2+y^2} 
  \right) .
\]
Let $x_b>0$ be such that $\frac{-1+x_b^2}{1+x_b^2} = \cos \theta_b$, i.e.\ $x_b = \tan ((\pi - \theta_b)/2)=1/q$.
Then we choose 
 \[
   \overline{m}_0(r;s) = T \left( x_b + x_b\cos( 2\pi s) \frac{1-r}{r}  , 
  x_b\sin( 2\pi s) \frac{1-r}{r} \right),
 \]
see also Figure~\ref{f:cover1}. 
For the special case that $m_b=N$ we choose
\begin{equation}
\label{e:mN}
  \overline{m}_0^N(r;s) = T\bigl( \tan(\pi (r-1/2)),\tan(\pi(s-1/2)) \bigr).
\end{equation}
This is just one explicit choice; any homotopy of this family of initial data that obeys the boundary conditions (\ref{e:endpoints}) also represents a degree 1 map on $S^2$.
Let $X_0$ be the collection of initial data obtained by taking all such homotopies. It is not hard to see that $X_0$ is the space of continuous functions (with the usual supremum norm) satisfying the boundary conditions.
Let $X_1$ be the subset of initial data in $X_0$ for which the solution to the equivariant equation~(\ref{LLG_3comp}) blows up in finite or infinite time.
The following result states that the co-dimension of $X_1$ is at most one. In particular, each one parameter family of initial data that represent a \emph{degree 1} map from $S^2$ to itself has at least one member that blows up. 

\begin{proposition}
\label{prop:one}
  Let $\beta>0$.
  The blowup set $X_1$ for the equivariant equation~(\ref{LLG_3comp}) has co-dimension at most 1.
\end{proposition}
\begin{proof}
Let $m_0(r,s)$ be \emph{any} family of initial data that, via the above identification, represent a degree 1 map from $S^2$ to itself.
Let $m(t,r;s)$ correspond to the solution with initial data~$m_0(r;s)$.
As explained above, we see from (\ref{e:endpoints}) that we may view $m_0(\cdot;\cdot)$ as a map from $S^2 \to S^2$ by identifying $\{0\} \times S^1$ and $\{1\} \times S^1$ to points.
Since the boundary points are fixed in time, we may by the same argument
view $m(t,\cdot;\cdot)$ as a map from $S^2$ to itself along the entire evolution.
Note that since the energy is a Lyapunov functional for $\beta>0$, any solution tends to an equilibrium as $t \to \infty$.
If none of the solutions in the family
would blow up (in finite or infinite time) along the evolution, then all
solutions converge smoothly to the unique equilibrium. 
In particular, for large $t$ the map $m(t,\cdot;\cdot): S^2 \to S^2$ has its image in a small neighborhood of this equilibrium, hence it is contractible and
thus has degree 0. Moreover, if there is no blowup then $m(t,\cdot;\cdot) : S^2 \to S^2$ is continuous in $t$,  i.e.\ a homotopy. 
This is clearly contradictory, and
we conclude blowup must occur for at least one solution in the
one-parameter family $m_0(r,s)$. This is a topologically robust property in the sense that any small perturbation of $m_0(r,s)$ also represents a degree 1 map, and the preceding arguments thus apply to such small perturbations of    $m_0(r,s)$ as well. This proves that the co-dimension of $X_1$ is at most 1. 
\end{proof}

One may wonder what happens when (equivariant) symmetry is lost. Although a~priori 
the co-dimension could be higher in that case, we will show that this is
not so.
For convenience, 
we only deal with boundary conditions $m(t,x) = N$ for all $x \in \partial D$, which simplifies the geometric picture,
but the argument can be extended to more general boundary conditions.

The only equilibrium solution in this
situation is $m(x)\equiv N$~\cite{Lemaire}.
Let $\overline{m}_0^N(r;s)$ be the family of initial data for the equivariant case with
boundary condition $m_b=N$ (see (\ref{e:mN}) and Figure~\ref{f:cover1}).
Consider now the following family of initial data for the general case:
\[
  \overline{M}_0^N(x;s)=\overline{M}_0^N(r,\psi;s)= 
  \left( \begin{array}{ccc}
  \cos \psi & -\sin \psi & 0\\
  \sin \psi &  \cos \psi & 0\\
  0 & 0 & 1 
   \end{array}\right) 
  \overline{m}_0^N(r;s) .
\]
We see that  $\overline{M}_0^N$ maps $D \times [0,1]$ to $S^2$, and $\overline{M}_0^N(\partial D; [0,1])=N$,
but also $\overline{M}_0^N(D^2;0)=\overline{M}_0^N(D^2;1)=N$. 
We may thus identify $\partial D \times [0,1] \cup D
\times \{0,1\}$ to a point, and interpret $\overline{M}_0^N$ as a map from $S^3$ to $S^2$.
In particular, $\overline{M}_0^N$ represents an element in the homotopy group
$\pi_3(S^2) \cong \mathbb{Z}$. Furthermore, upon inspection, $\overline{M}_0^N$
represents the generator of the group, since it is (a deformation of) 
the Hopf map (see e.g.~\cite{Hatcher}). 
Let $\widetilde{X}_0$ be the collection of initial data in one parameter families obtained from all homotopies of $\overline{M}_0^N$ that obey the boundary conditions
\begin{equation}
\label{e:cylinder}
  \overline{M}_0(\partial D; [0,1])=N, \qquad\text{and}\qquad 
  \overline{M}_0(D^2;0)=\overline{M}_0(D^2;1)=N. 
\end{equation}
Let $\widetilde{X}_1$ be the subset of initial data in $\widetilde{X}_0$ for which the solution to the differential equation~(\ref{LLG1}) blows up in finite or infinite time.
As before, the co-dimension of $\widetilde{X}_1$ is at most 1, showing that dropping the equivariant symmetry does not further increase the instability of the blowup scenario.

\begin{proposition}
	\label{prop:two}
  Let $\beta>0$ and $m_b=N$.
  The blowup set $\widetilde{X}_1$ for the general equation~(\ref{LLG1}) has co-dimension at most 1.
\end{proposition}
\begin{proof}
The proof is analogous to the one of Proposition~\ref{prop:one}, but one uses $\pi_3(S^2)$ instead of $\pi_2(S^2)$, i.e.\ the degree, to obtain the contradiction.
\end{proof}

As a final remark, even though blowup is co-dimension one, this does not
mean it is irrelevant. Clearly, by changing the initial data slightly one
may avoid blowup. On the other hand, the arguments above indicate that 
blowup is caused by the topology of the target manifold, 
and one can therefore not circumvent this type of
singularity formation by simply adding additional terms to the equation
(for example a physical effect that works on a smaller length scale),
unless additional equilibria are introduced which reflect the pinning of a defect.

\section{Asymptotic analysis}
\label{sec:Asymp}
In this section we extend the results of \cite{BHK}, where the rate
of blowup for radially symmetric solutions to the harmonic heat map
problem (\ref{radHM}) was determined.
  We will consider both the extension to the full
Landau-Lifshitz-Gilbert equation (i.e.~$\alpha \ne 0$), as well as allowing a
particular class of non-radial perturbations.  We find that
blowup solutions are always unstable in the equivariant regime.  It can be
understood that the blowup solutions are separatrices between two distinct
global behaviors.

\subsection{The inner region}
We will proceed with an expansion motivated by two facts: (i) blowup in
the harmonic map heatflow is a quasi-static modulated stationary solution;
(ii) the full LLG problem has the same stationary profiles as the harmonic
map heatflow.

Without specifying the rescaling factor $R(t)$ yet we introduce the 
rescaled variable
\[ \xi = \frac{r}{R(t)}\]
which transforms (\ref{LLG_2comp}) to
\begin{equation}
\begin{aligned}
\theta_{\xi \xi} + \frac{1}{\xi} \theta_{\xi} - \frac{\sin 2\theta}{2}
\left(\frac{1}{\xi^2} + \phi_\xi^2\right)  &= \beta \left(R^2 \theta_t - R R' \xi \theta_\xi \right) + 
\alpha \sin \theta \left(R^2 \phi_t - R R' \xi \phi_\xi \right) , \\
 \phi_{\xi \xi} + \frac{1}{\xi} \phi_\xi + \frac{\sin 2\theta}{\sin^2 \theta} \phi_\xi \theta_\xi &= 
 \beta \left (R^2 \phi _t - R R' \xi \phi_\xi \right) - \frac{\alpha}{\sin \theta} 
 \left(R^2 \theta_t - R R' \xi \theta_\xi \right) .
\end{aligned}
\label{llgEq2}
\end{equation}
To solve this in the limit $R \to 0$ we pose the expansion
\begin{align*}
\theta & \sim \theta_0 + (\beta R R' - \alpha R^2 C') \theta_1 + \ldots , \\
\phi & \sim \phi_0 + (\beta R^2 C' + \alpha R R')\phi_1 + \ldots ,
\end{align*}
where
\begin{align}
\theta_0 & = 2 \arctan \xi ,\\
\phi_0 & = C(t),
\end{align}
represent slow movement along the two-parameter family of equilibria $\theta=2\arctan({r}/{R})$ and $\phi=C$. 


At the next order we have
\begin{align}
& \frac{d^2 \theta_1}{d\xi^2} + \frac{1}{\xi}\frac{d\theta_1}{d\xi} - \frac{\cos 2\theta_0}{\xi^2}\theta_1
= -\xi \frac{d \theta_0}{d\xi} , \\
&\frac{d^2 \phi_1}{d\xi^2} + \frac{1}{\xi}\frac{d\phi_1}{d\xi} + \frac{\sin 2\theta_0}{\sin^2\theta_0} \frac{d \theta_0}{d\xi} \frac{d \phi_1}{d\xi} = 1 .
\label{e:phi1}
\end{align}
These equations can both be solved exactly, but we omit the algebraic details since for the matching we only need the asymptotic behaviour as $\xi \to \infty$, viz.:
\begin{alignat}{2}
\theta_1 &\sim -\xi \ln \xi + \xi,&& \qquad \mbox{as } \xi \to \infty,\\
\phi_1 &\sim \frac{1}{2} (\ln \xi) \xi^2 - \frac{1}{2}\xi^2, &&\qquad \mbox{as } \xi \to \infty.
\end{alignat}

\begin{figure}[t]
{\includegraphics[width=\textwidth]{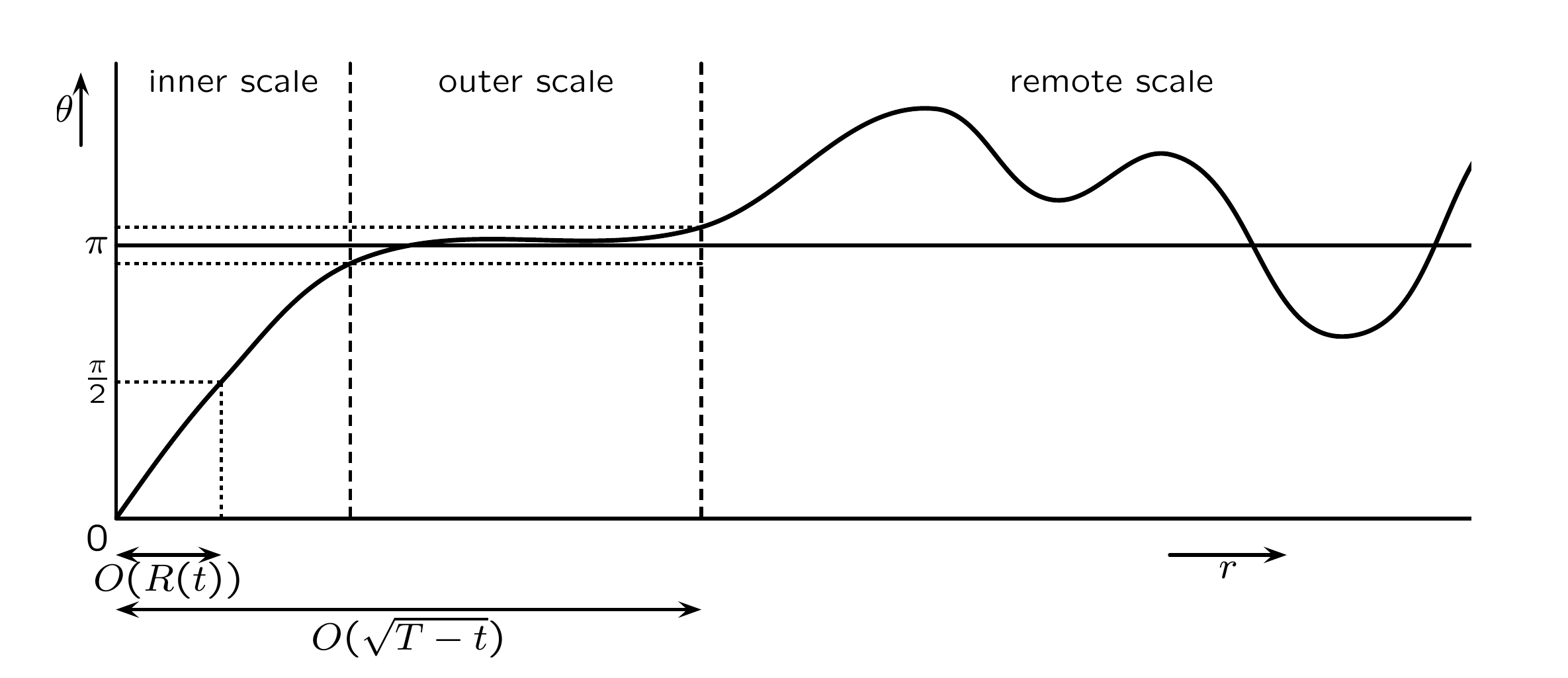}}
\caption{Regions for asymptotic analysis}
\label{fig:regions}
\end{figure}

At this stage both $R(t)$ and $C(t)$ are unspecified functions.  They will be determined through the matching of the inner ($r\sim R(t)$) and outer regions ($r \sim \sqrt{T-t}$), cf.\ Figure~\ref{fig:regions}.

\subsection{The outer region}
\label{sec:outer}

To make the mechanics of the linearization and matching as transparent as
possible we shall now change variables by linearizing around the south pole in
the formulation (\ref{LLG_3comp}):
\[ (\pi - \theta) e^{i\phi} = u+i v,  \qquad w = -1.\]
This recovers
\begin{align}
u_t & = \beta \left(u_{rr} + \frac{1}{r}u_r - \frac{1}{r^2} u \right) + \alpha \left(v_{rr} + \frac{1}{r}v_r - 
\frac{1}{r^2} v \right) \\
v_t & = \beta \left(v_{rr} + \frac{1}{r}v_r - \frac{1}{r^2} v \right) - \alpha \left(u_{rr} + \frac{1}{r}u_r - 
\frac{1}{r^2} u \right)
\end{align}
To solve this we introduce
$ z = u + i v$,
whence 
\begin{equation}\label{Zeq1}
z_t = (\beta - i \alpha) \left(z_{rr} + \frac{1}{r}z_r - \frac{1}{r^2} z\right).
\end{equation}
This is simply the projection of the flow from the sphere on to the tangent
plane at the pole.  Notice that in the respective limits we the recover
modified linear heat equation ($\alpha=0$) and Schr\"odinger equation ($\beta=0)$ on the tangent plane as appropriate.

To match the inner and outer regions we first define the similarity variables
\begin{equation}
\label{SimVars}
\xi = \frac{r}{R(t)}, \quad \tau = -\ln(T-t), \quad y = e^{\tau/2} r.
\end{equation}
To reduce confusion in what follows we shall denote 
\[ \frac{d f(t)}{dt} = f' \qquad \mbox{and} \qquad \frac{d f(\tau)}{d \tau} = \dot{f}, \qquad \mbox{thus} \quad
\frac{d}{dt} f(\tau) = e^\tau \dot{f}. \]

Under this change of variables equation (\ref{Zeq1}) becomes
\begin{equation}\label{ZeqSimVars}
z_\tau = {\mathcal L} z \equiv -\frac{y}{2} z_y +  (\beta - i \alpha) \left(z_{yy} + \frac{1}{y}z_y - \frac{1}{y^2} z\right).
\end{equation}


A solution for this equation is
\[ z =  \sigma e^{-\tau/2} y \]
for any $\sigma$ --- this is just $z = \sigma r$ in (\ref{Zeq1}).  This corresponds to
the eigenfunction of the dominant eigenvalue of ${\mathcal L}$, which governs the generic
long time behaviour of solutions of~(\ref{ZeqSimVars}).   
When we allow $\sigma$ to vary slowly with $\tau$, we obtain a series expansion
for the solution of the form 
\begin{equation}\label{ZeqMatch}
z \sim e^{-\tau/2}\left[ \sigma(\tau) y + \dot{\sigma}(\tau)\left((\beta -
    i \alpha)\frac{4}{y} - 2 y \ln y \right) + \ldots \right] 
\qquad \mbox{as} \; \tau \to \infty .
\end{equation}
We now see that the introduction of $\alpha$ non-zero does not affect the
procedure for the expansion. 

Denoting $\sigma(\tau) = \sr(\tau) + i \si(\tau)$, we introduce
\begin{equation}
	\begin{aligned}
\lr + i \li & \equiv 
(\dsr\beta + \dsi \alpha) + i (\dsi \beta - \dsr \alpha)   \\ &=
(\beta - i \alpha) (\dsr + i \dsi) .
\end{aligned}
\label{LamDef}
\end{equation}
We recover $\theta$ and $\phi$ through $|z| = \pi - \theta$, $\arg z = \phi$ and expand for
small $y$:
\begin{alignat}{2}
\nonumber
|z| & = e^{\tau/2} \sqrt{ \left(\sr y + \lr \frac{4}{y}+ \dots\right)^2  + 
\left(\si y + \li \frac{4}{y}+ \dots\right)^2 } , \\
 \pi-\theta & \sim e^{-\tau/2} \left(\sqrt{\lr^2 + \li^2}\frac{4}{y} +
   \frac{\sr\lr + \si\li}{\sqrt{\lr^2 + \li^2}} y +\dots \right) 
&&\text{for small } y,
 \label{ThetaOutAsy} \\
 \nonumber 
\arg z & = \arctan\left(\frac{\si y + 4\li y^{-1}+ \dots}{\sr y + 4 \lr
    y^{-1}+ \dots}\right) , \\ 
\phi & \sim \arctan{\frac{\li}{\lr}} + 
\frac{\si\lr-\sr\li}{\lr^2 + \li^2} \frac{y^2}{4} + \dots 
&&\text{for small } y.
\label{PhiOutAsy}
\end{alignat}
Here and in what follows one should be slightly careful interpreting all formulae involving the $\arctan$ because of multi-valuedness. 
For future reference we note that
\begin{alignat}{2}
  \arg z &\to  \arctan\frac{\si}{\sr} &\qquad&\text{for large } y, \label{e:angle1}\\
  \arg z &\to \arctan\frac{\li}{\lr} = 
         \arctan\frac{\dsi}{\dsr} - \arctan\frac{\alpha}{\beta}
   &&\text{for small } y. \label{e:angle2}
\end{alignat}

\subsection{The matching}

In order to match the inner region to the outer we first write the inner solution
in the similarity variables: 
\begin{align}
\nonumber
\theta & \sim 2 \arctan\left(e^{-\tau/2}\frac{y}{R}\right)  \\
\nonumber
& \hspace*{1cm} + e^\tau\left(\beta R \dot{R} - 
\alpha R^2 \dot{C}\right)\left(-e^{-\tau/2} \frac{y}{R} \ln\left(e^{-\tau/2}\frac{y}{R}\right) +
e^{-\tau/2}\frac{y}{R}\right) + \ldots \\
\label{ThetaInAsy}
& \sim \pi - 2 \frac{R e^{\tau/2}}{y} + e^{\tau/2} \left(\beta \dot{R} - \alpha R \dot{C}\right)
\left(\frac{\tau}{2} + \ln R - 1\right) y  + \ldots , \\
\nonumber
\phi & \sim C + \left(\beta R^2 \dot{C} + \alpha R \dot{R}\right) e^\tau\left(\frac{e^{-\tau} y^2}{R^2}
\ln\left(\frac{e^{-\tau/2}y}{R}\right) - \frac{y^2 e^{-\tau}}{R^2}\right)  + \ldots \\
\label{PhiInAsy}
&\sim C - \frac{1}{2}\left(\beta \dot{C} + \alpha \frac{\dot{R}}{R}\right) 
\left( \frac{\tau}{2} + \ln R - 1\right)y^2  + \dots .
\end{align}
The matching procedure now involves setting $C$ and $R$ such that the
expansions (\ref{ThetaInAsy}) and (\ref{PhiInAsy}) agree with
(\ref{ThetaOutAsy}) and (\ref{PhiOutAsy}) respectively, to two orders in
$y$:
\begin{alignat*}{3}
\theta&:\quad  &\cO(y^{-1}) &: ~ & 2 R e^{\tau/2} &\sim 4 e^{-\tau/2}
(\lr^2 + \li^2)^{1/2} , \\
&&\cO(y) &:  & 
-\left(\beta \dot{R} - \alpha R \dot{C}\right) 
\left( \frac{\tau}{2} + \ln R \right) e^{\tau/2}
&\sim e^{-\tau/2} \frac{\sr\lr + \si\li}{(\lr^2 + \li^2)^{1/2}} , \\
\phi&:&  \cO(y^0)&:  &  C &\sim \arctan\left(\frac{\li}{\lr}\right) ,
\\
&& \cO(y^2)&: &  
-\frac{1}{2}\left(\beta \dot{C} + \alpha  \frac{\dot{R}}{R}\right)
\left(\frac{\tau}{2} + \ln R  \right) 
& \sim \frac{1}{4} \frac{ \si\lr -  \sr\li}{\lr^2 + \li^2} .
\end{alignat*}
To solve this we set $R = e^{-\tau} p(\tau)$ (with $p(\tau)$ algebraic in $\tau$), and after rearranging terms we
get
\begin{equation}
	\left\{
\begin{aligned}
  & p \sim 2  (\lr^2 + \li^2)^{1/2} , \\
  & C \sim \arctan\left(\frac{\li}{\lr}\right) , \\
 & \left(\beta (\dot{p}-p) - \alpha p \dot{C}\right) 
\left( \frac{\tau}{2} - \ln p \right)
\sim \frac{\sr\lr + \si\li}{(\lr^2 + \li^2)^{1/2}} , \\
& \left(\beta p \dot{C} + \alpha  (\dot{p}-p) \right) 
\left(\frac{\tau}{2} - \ln p  \right) 
 \sim \frac{ \si\lr - \sr\li}{(\lr^2 + \li^2)^{1/2}} .
\end{aligned}
\right .
\label{FullSys1}
\end{equation}
Since $p$ is defined not to change exponentially fast in $\tau$, we
neglect the terms of $\cO ( \ln p)$.
Using the definition of $\lambda$ we may simplify (\ref{FullSys1})  to get
\begin{equation}
	\left\{
\begin{aligned}
  p &\sim 2  (\dsr^2 + \dsi^2)^{1/2} , \\
  C &\sim \arctan\frac{\dsi}{\dsr}-\arctan\frac{\alpha}{\beta} , \\
 \frac{\tau}{4}p (\dot{p}-p) 
&\sim \beta(\sr\lr + \si\li)+\alpha(\si\lr - \sr\li)  = \sr\dsr+\si\dsi \,, \\
\frac{\tau}{4}p^2\dot{C}
& \sim \beta(\si\lr - \sr\li) -\alpha(\sr\lr + \si\li)= \si\dsr - \sr\dsi \,.
\end{aligned}
\right .
\label{FullSys2}
\end{equation}
Finally, we introduce $\tC=C+\arctan\frac{\alpha}{\beta}$, so that
\begin{equation}
	\left\{
\begin{aligned}
  p &\sim 2  (\dsr^2 + \dsi^2)^{1/2} , \\
  \tC &\sim \arctan\frac{\dsi}{\dsr} \,, \\
 \frac{\tau}{4}p (\dot{p}-p) 
& \sim \sr\dsr+\si\dsi \,, \\
\frac{\tau}{4}p^2\dot{\tC} 
& \sim  \si\dsr - \sr\dsi \,,
\end{aligned}
\right .
\label{FullSys3}
\end{equation}
which is \emph{independent} of $\alpha$ and $\beta$. This formulation strongly suggests
that the case $\alpha = 1, \beta = 0$ is not different from the dissipative case $\beta > 0$. Before solving and studying the system (\ref{FullSys3}) let us recall what its solutions tell us: $p(\tau)$ gives an algebraic correction to the blowup rate, $\tC(\tau)$ determines the local behaviour of $\phi$ near blowup and $\sigma$ describes the amplitude and orientation of the solution in self-similar coordinates, see (\ref{e:angle1}),(\ref{e:angle2}).  In order to fully understand blowup, we need to solve for the blowup coordinates and determine their stability.

\begin{figure}
\centerline{\includegraphics[height=5cm]{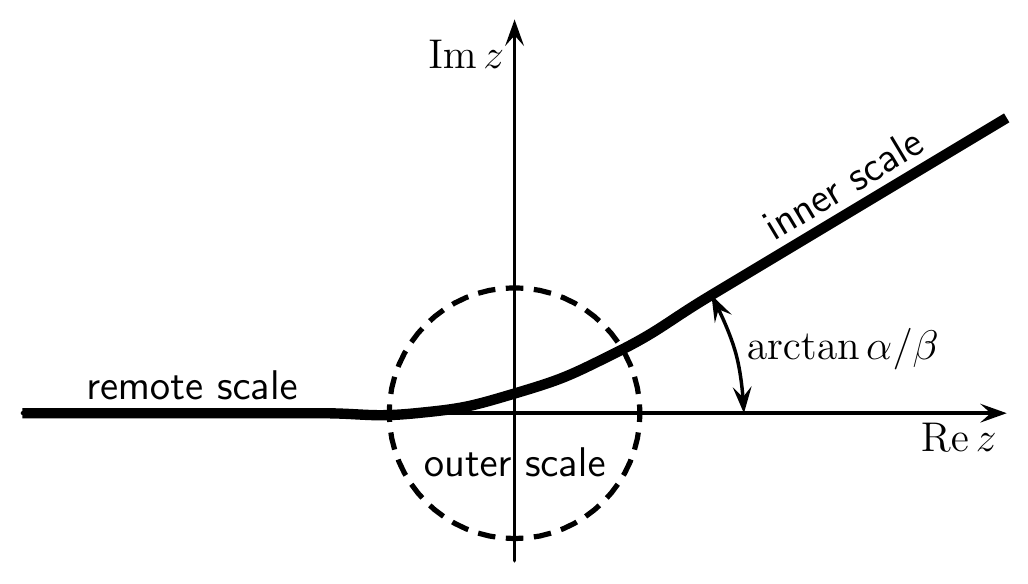}}
\caption{The tangent plane at the south pole can be identified with the complex plane. The thick curve represents $z(t,\cdot)$ for a time near blowup. The angle between the solution in the inner scale (which bubbles off) and the remote scale (which remains) is indicated.}	
\label{fig:angle}
\end{figure}

The blowup solution is represented by $\si(\tau) = c \sr(\tau)$ for some
constant $c\in\mathbb{R}$ (or $c=\infty$, i.e.~$\sr=0$), with $\tan \tC= c$. In particular,
equations~(\ref{e:angle1}),(\ref{e:angle2}) show that there is an angle
\mbox{$\pi - \arctan \frac{\alpha}{\beta}$} between the sphere bubbling off and the
solution remaining at/after blowup, see Figure~\ref{fig:angle}.

By rotating the sphere we may take $\tC=0$ without loss of generality, i.e.\ $\si=\dsi=0$ and 
$\dsr >0$ (note the sign), see (\ref{FullSys3}). 
The blowup dynamics is described by 
\begin{equation}
\label{e:sigma0}
(\ddot{\sigma}_r - \dot{\sigma}_r) \tau = \sigma_r \,,
\end{equation}
and $p= 2\dsr >0$.
This equation has general solutions of the form
\[ 
\sr = k_1 \tau e^{\tau} + k_2 f(\tau) \qquad \mbox{where~} f(\tau)
\sim \frac{1}{\tau} \mbox{ as }  \tau \to \infty.
\] 
We can immediately
set $k_1 \equiv 0$ as this ``instability'' reflects shifts in the blowup
time and hence is not a real instability --- this is common to all blowup
problems~\cite{VAGbook}. 
We note that $k_2<0$ so that indeed $\dot{\sigma}_r > 0$ as $\tau\to\infty$,
and $p \sim -2 k_2 \tau^2 >0$.
This implies that $\sigma_r < 0$, hence $\arg z \to \pi$ for large $y$ (cf.\ Figure~\ref{fig:angle}). 

At this stage we have an asymptotic description of the blowup rate and its local
structure.  Unfortunately, we do not have enough information to determine stability.
To more carefully understand the dynamics of this system we need
to linearize about this leading order solution to find the subsequent
corrections $\sr[1]$, $\si[1]$, $p_1$ and $\tC_1$ in $\sr$, $\si$, $p$
and $\tC$, respectively. Taking 
$\sr[0]= f(\tau) \sim k_2/\tau$, $\si[0]=0$,
$p_0=2f'(\tau) \sim -2 k_2 /\tau^2$, $\tC_0 =0$, 
we get as the system for the  next order
\begin{equation}
	\left\{
\begin{aligned}
  p_0 p_1 & = 4 \dsr[0]\dsr[1] \,, \\
  \tC_1 &= \frac{\dsi[1]}{\dsr[0]} \,, \\
 \frac{\tau}{4} (p_1 \dot{p}_0 + p_0 \dot{p}_1- 2p_0 p_1) 
& = \sr[0]\dsr[1]+\sr[1]\dsr[0] \,, \\
\frac{\tau}{4}p_0^2\dot{\tC}_1 
& =  \si[1]\dsr[0] - \sr[0]\dsi[1] \,,
\end{aligned}
\right .
\label{FullSecOrd}
\end{equation}
which separates into two systems. The first one is (using $p_0=2\dsr[0]$)
\begin{equation}
	\left\{
\begin{aligned}
   p_1 & =  2 \dsr[1] \,, \\
   \frac{\tau}{2} (p_1 \ddot{\sigma}_{0r} + \dsr[0] \dot{p}_1- 2\dsr[0] p_1) 
& = \sr[0] \dsr[1]+ \dsr[0] \sr[1] \,,
\end{aligned}
\right .
\label{SecOrd1}
\end{equation}
which, using that $\sr[0]$ solves (\ref{e:sigma0}), reduces to
\[
(\ddot{\sigma}_{1r} - \dsr[1]) \tau = \dsr[1] \,,
\]
the same equation as for $\sr[0]$, and provides no additional information.
The other system is
\begin{equation}
	\left\{
\begin{aligned}
   \tC_1 &= \frac{\dsi[1]}{\dsr[0]} , \\
 \tau \dsr[0]^2 \dot{\tC}_1 
& =  \si[1]\dsr[0] - \sr[0]\dsi[1] \, ,
\end{aligned}
\right .
\label{SecOrd2}
\end{equation}
which can be rewritten as
\[
\tau(\dsr[0]\ddot{\sigma}_{1i} - \dsi[1]\ddot{\sigma}_{0r} ) =
\si[1]\dsr[0] - \sr[0]\dsi[1] \,, 
\]
or, again using that $\sr[0]$ solves (\ref{e:sigma0}),
\[
(\ddot{\sigma}_{1i} - \dsi[1]) \tau = \dsi[1] \,,
\]
i.e., once again equation~(\ref{e:sigma0}). The asymptotic behaviour
of $\si[1]$ and $\tC$ is thus given by ($\kappa_1,\kappa_2 \in \mathbb{R}$)
\begin{alignat*}{1}
  \si[1]&\sim \frac{\kappa_1}{\tau}  + \kappa_2 \tau e^{\tau} , \\
  \tC_1& \sim - \tau \si[1] \sim - \kappa_1 -  \kappa_2 \tau^2 e^{\tau}, 
\end{alignat*}
where the exponentially growing terms show that blowup is unstable (the neutral mode corresponds to a (fixed, time-independent) rotation of the sphere).

\section{Numerical computations}
\label{sec:Num}


To supplement the formal analysis above we now present some numerical
experiments in the radial, equivariant and fully two-dimensional cases. 

\subsection{Numerical methods}
\label{sec:NumM}

To reliably numerically simulate potentially singular solutions to (\ref{LLG1}) one needs to use adaptivity in both time and space as well satisfy the constraint $|m(x,t)| = 1$. 
For the former, we use $r-$adaptive numerical methods as described in \cite{BuddWilliams}. This approach is based on the moving mesh
PDE approach of Huang and Russell \cite{HR} combined with scale-invariance
and the Sundman transformation in time.  The expository paper \cite{BuddWilliams} provides many examples of this method being effective for computing blowup solutions to many different problems. For the latter we can either use formulation (\ref{LLG_3comp}) and use a projection step or regularize the Euler angle formulation (\ref{LLG_2comp}).  We have
implemented both and found little difference in efficiency or accuracy and hence will use formulation (\ref{LLG_3comp}) for all but Example 1 as it directly
follows the above asymptotic analysis.

Full two-dimensional calculations have only been performed in the case of formulation (\ref{LLG1}) and on the unit disk.  This latter fact is for numerical convenience and in no way affects the structure of local singularities (should they arise). Here adaptivity was performed using the parabolic Monge-Ampere equation as described in
\cite{BW2}

\subsection{Numerical results}

\begin{figure}[!t]
{\includegraphics[width=\textwidth]{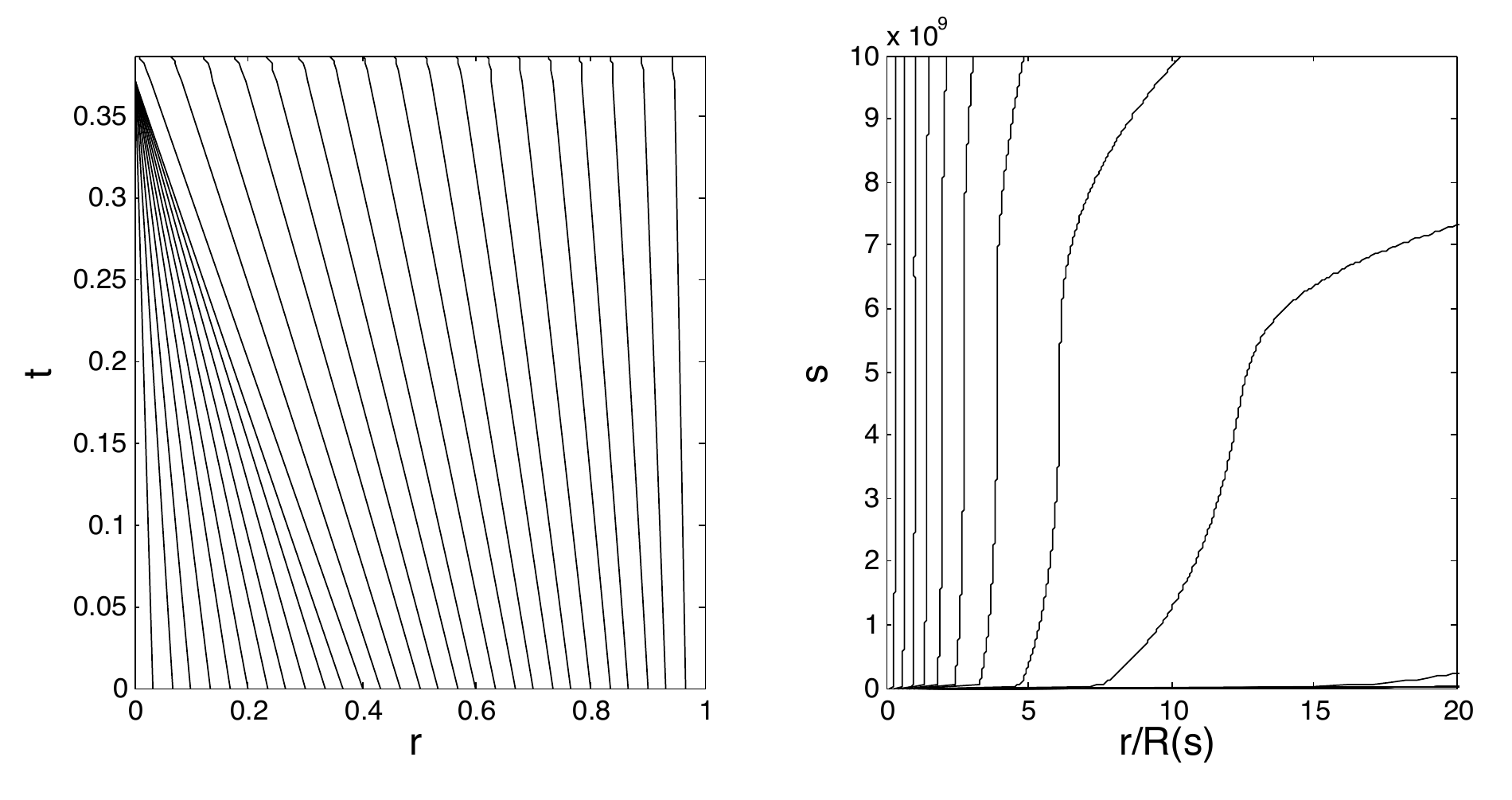}}
\caption{Left: Physical grid on which the solution was computed.  Right: Computational grid in the region of $\xi = 0$. Notice there is a region of essentially constant in $\xi$ grid trajectories in this region. Some trajectories are leaving this region as $\xi = r/R$ and $R \to 0^+$. Note, in both figures only every fifth grid trajectory is plotted.}
\label{fig:Example1a}
\end{figure}

\begin{figure}[!t]
{\includegraphics[width=\textwidth]{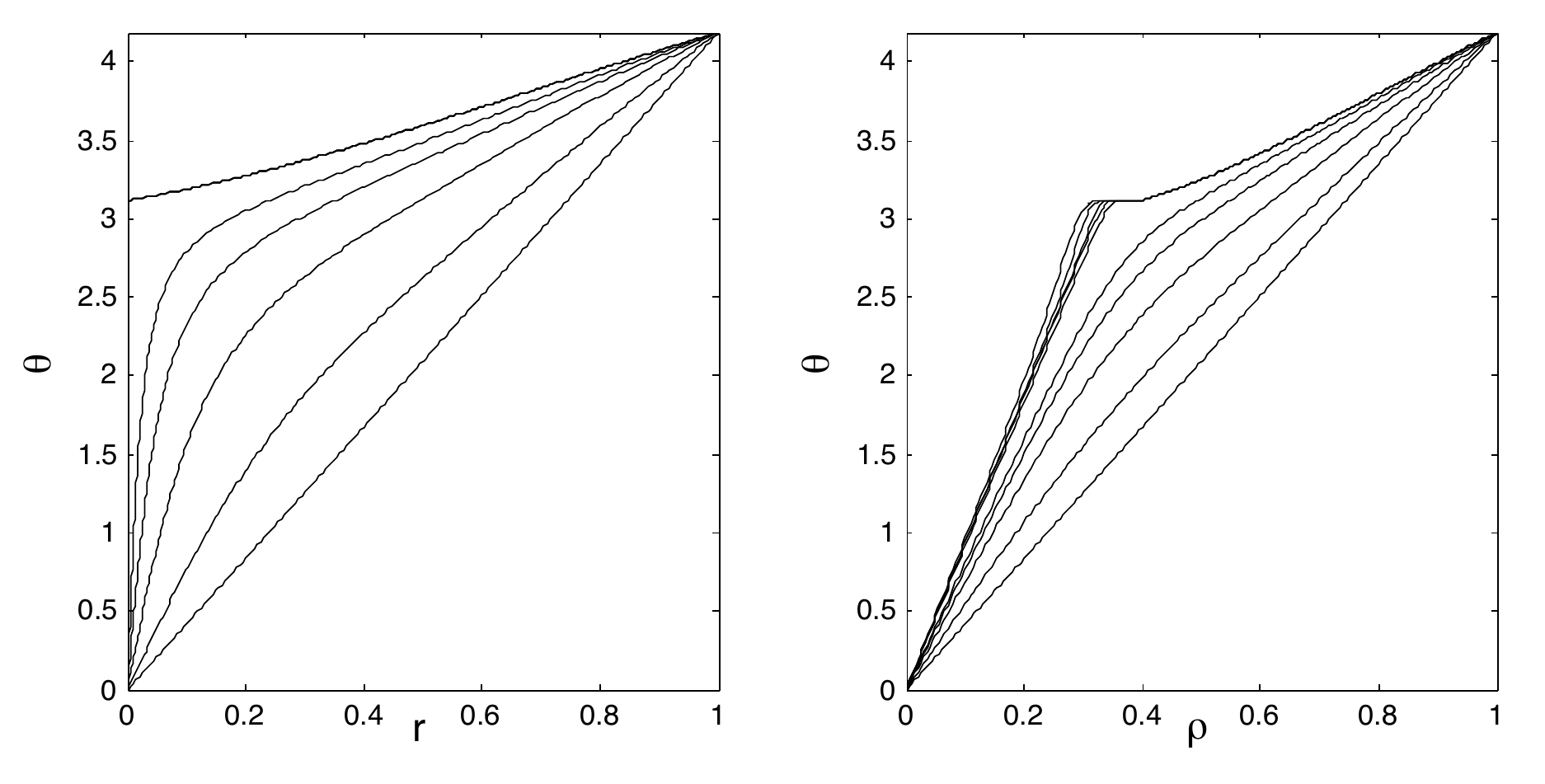}}
\caption{Left: Solution on physical grid at selected times. Right: Solution on
computational grid at same times. This clearly shows that the blowup region is very well resolved.}
\label{fig:Example1b}
\end{figure}

\label{sec:NumR}
In this Section we present a sequence of numerical experiments to validate the
results above. For examples 1-3 we used  $N=201$ spatial points, the monitor function
\[M = |\nabla m| + \int_{\Omega_c}|\nabla m| \, dx\]
and took
\[ \frac{dt}{ds} = \frac{1}{||M||_\infty}\]
as the rescaling between computational time $s$ and physical time $t$. Example 4 was computed on
a $61\times61$ grid using the same monitor functions. Here we took $|\nabla m| = \sqrt{u_r^2 + v_r^2 + w_r^2}$ when using form (\ref{LLG_3comp}) in radial coordinates,
 $|\nabla m| = |\theta_r|$ when using form (\ref{radHM}) and the Cartesian gradient when solving the problem in two
dimensions. In one dimension $\Omega_c = [0,1]$ and in two dimensions $\Omega_c$ is the unit disk. In both cases the integral
in the monitor function is computed in the physical variables.

{\bf Example 1 - the Radial harmonic map}
\begin{figure}[!t]
\center{\includegraphics[width=.8\textwidth]{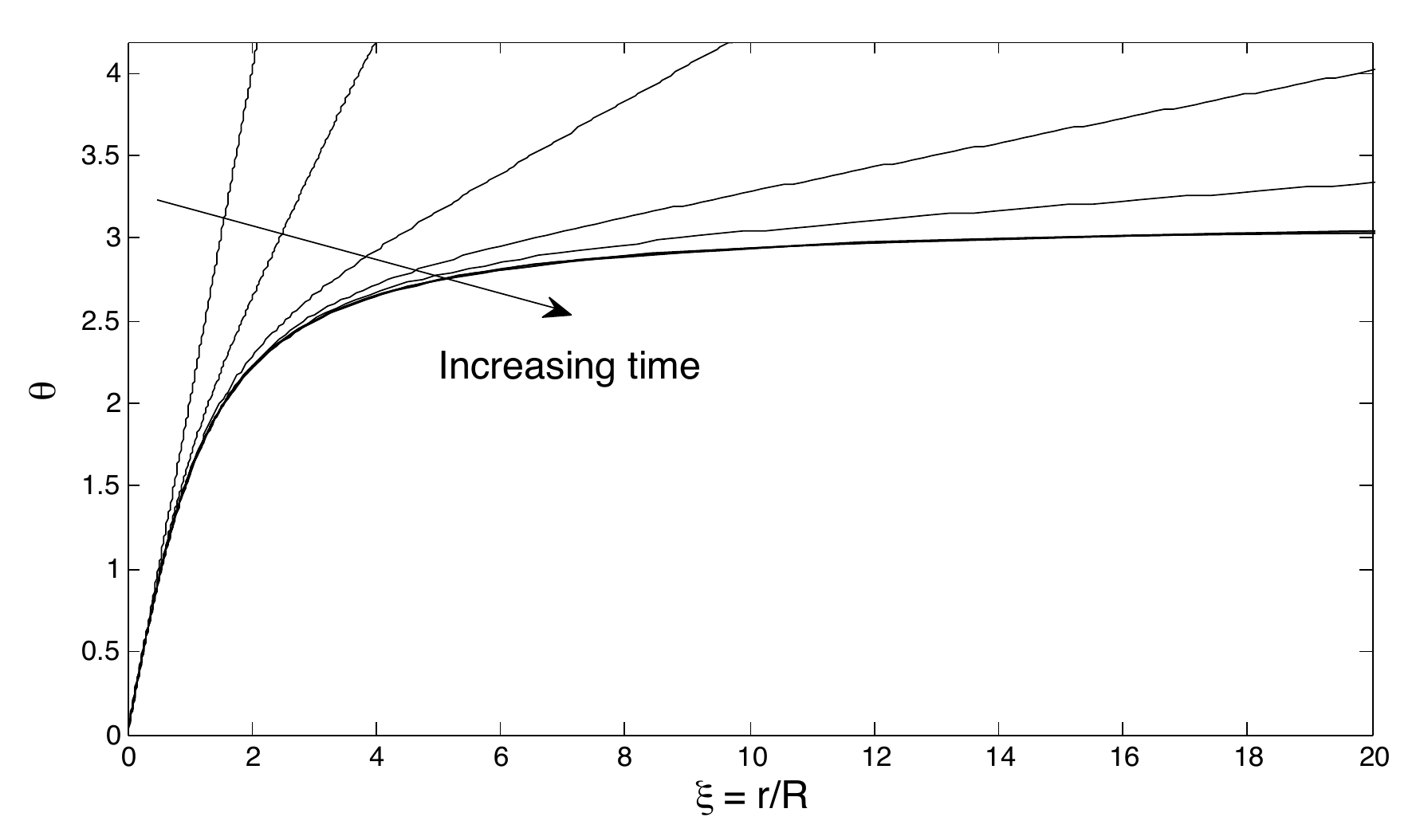}}
\caption{Evolution in the rescaled spatial variable. The solutions converge to the rescaled $\arctan$ profile, $\bar{\theta} = 2 \arctan (r/2R)$ as predicted (plotted under the numerical solutions).}
\label{fig:Example1c}
\end{figure}

The first example we will consider is equation (\ref{radHM}) with initial data
\begin{equation} 
	\theta_0(r) = \tfrac{4}{3} \pi r.
	\label{ThetaIC}
\end{equation}
	 This case has been proven to blowup with known structure \cite{CDY,Struwe} and asymptotically calculated rate $R(t)$. Figures \ref{fig:Example1a} and \ref{fig:Example1b} demonstrate the method and show how the
adaptive scheme follows the emerging similarity structure in the underlying evolution.
Figure \ref{fig:Example1c} shows excellent agreement with the analytical prediction of convergence to the $\arctan$ profile with $R(t)$ changing over twelve orders of magnitude.

{\bf Example 2 - Equivariant Harmonic map}
\begin{figure}[!t]
\centerline{{\includegraphics[width=.495\textwidth]{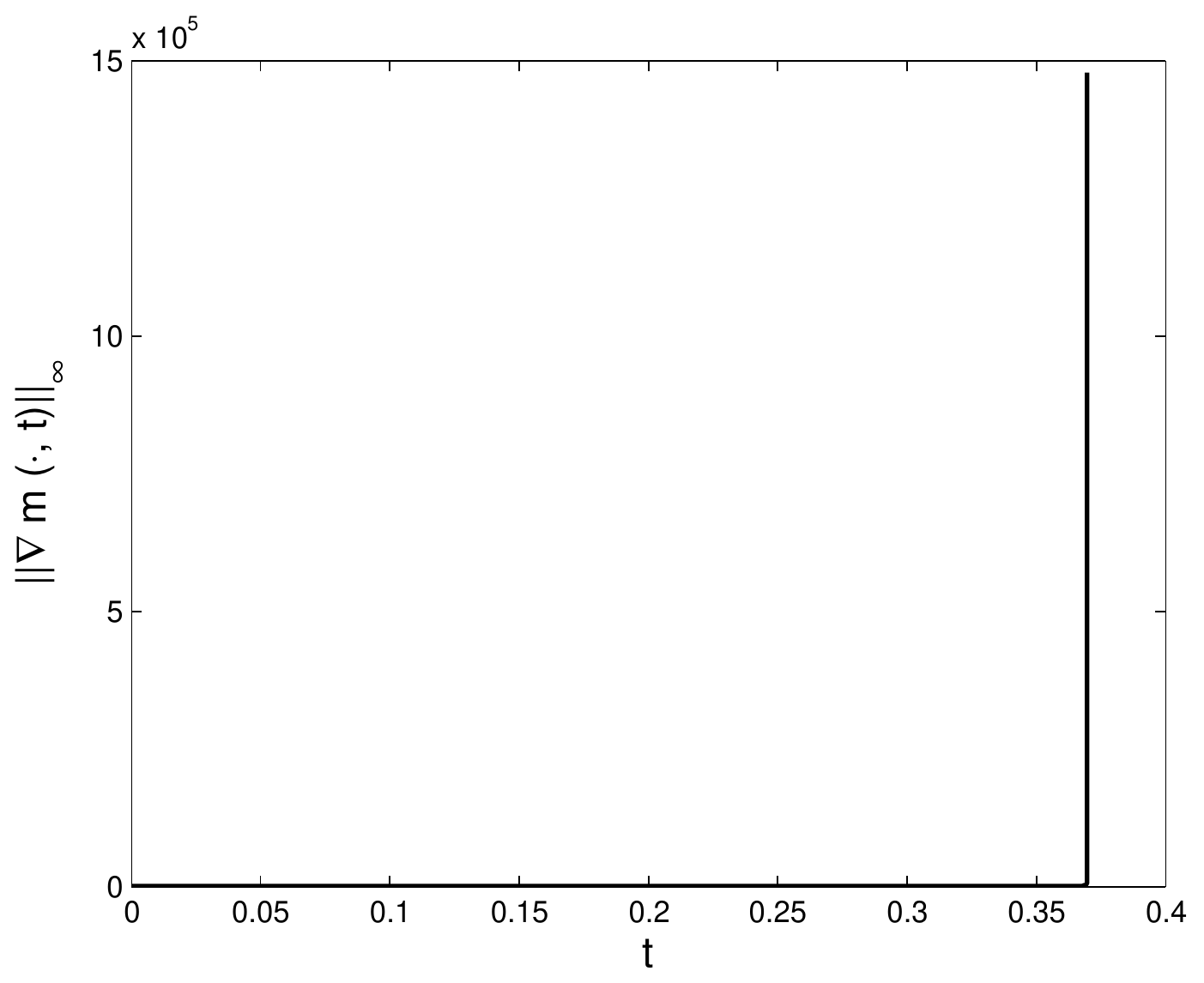}}
{\includegraphics[width=.5\textwidth]{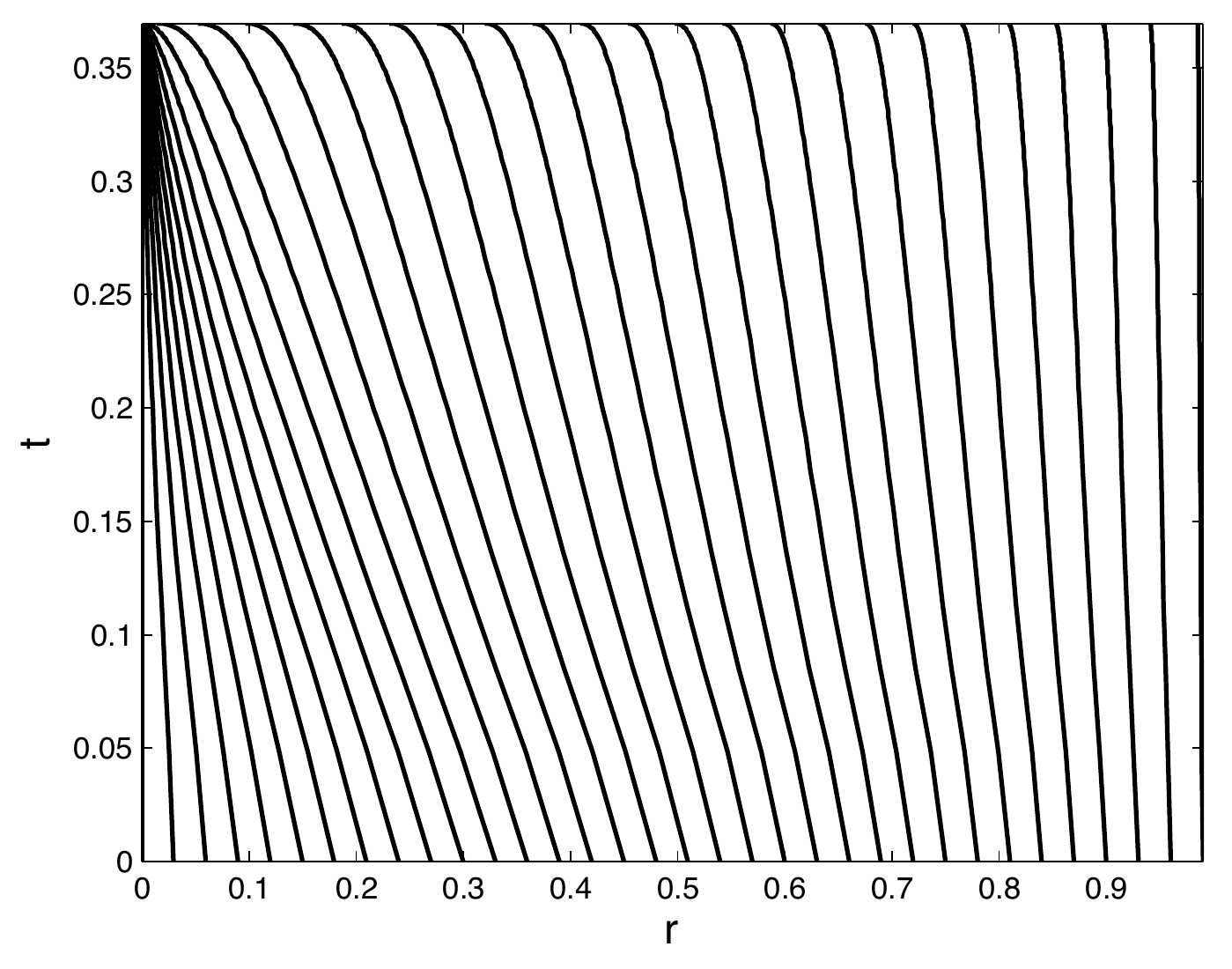}}}
\caption{Blowup of initial data (\ref{ThetaIC}) computed using (\ref{LLG_3comp}). Left: Evidence of blowup.
Right: Computational grid. Note that is very similar to Figure \ref{fig:Example1a} except that we cannot 
compute as far into the blowup.}
\label{fig:Example2a}
\end{figure}

\begin{figure}[!t]
\centerline{{\includegraphics[width=.485\textwidth]{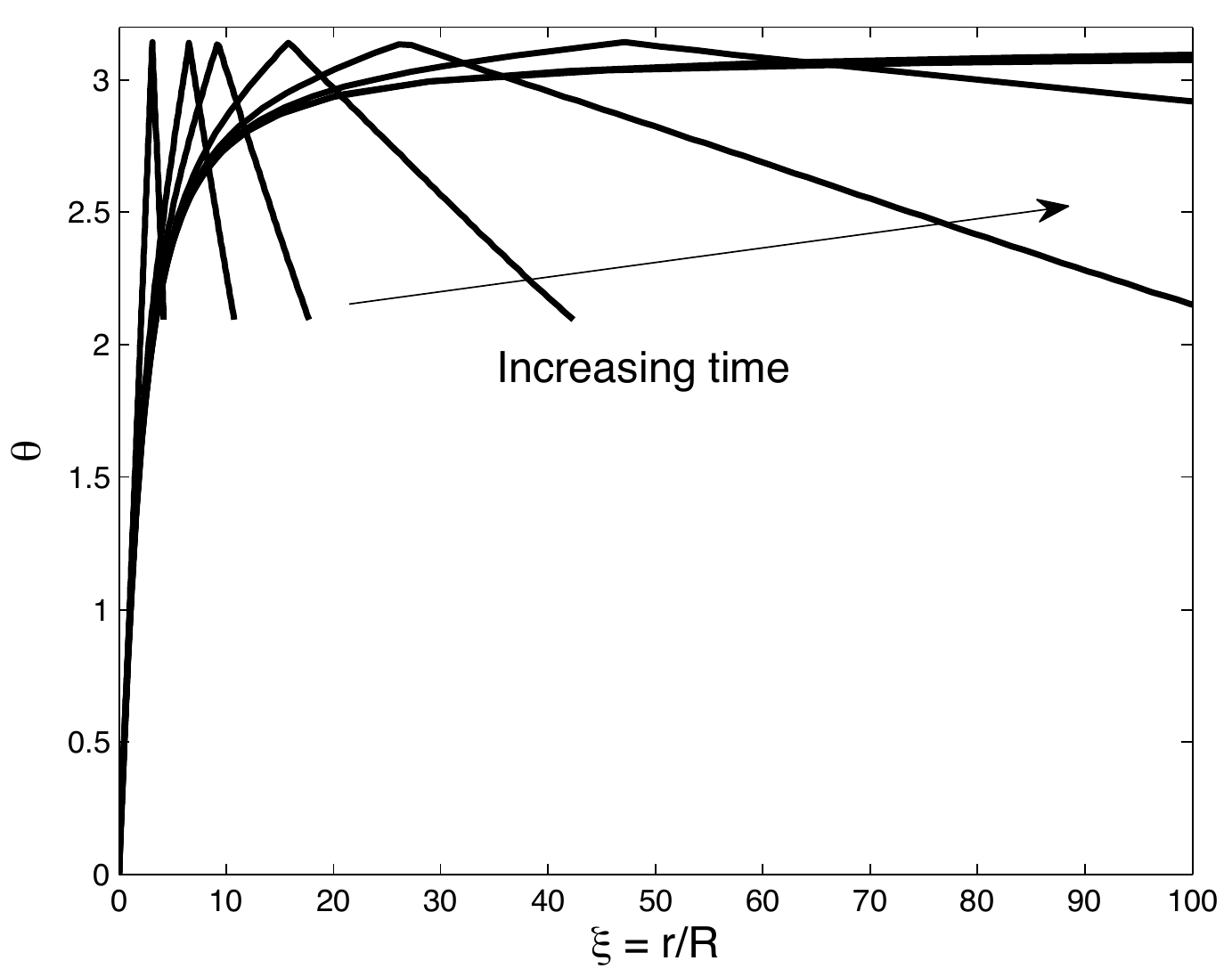}}
{\includegraphics[width=.49\textwidth]{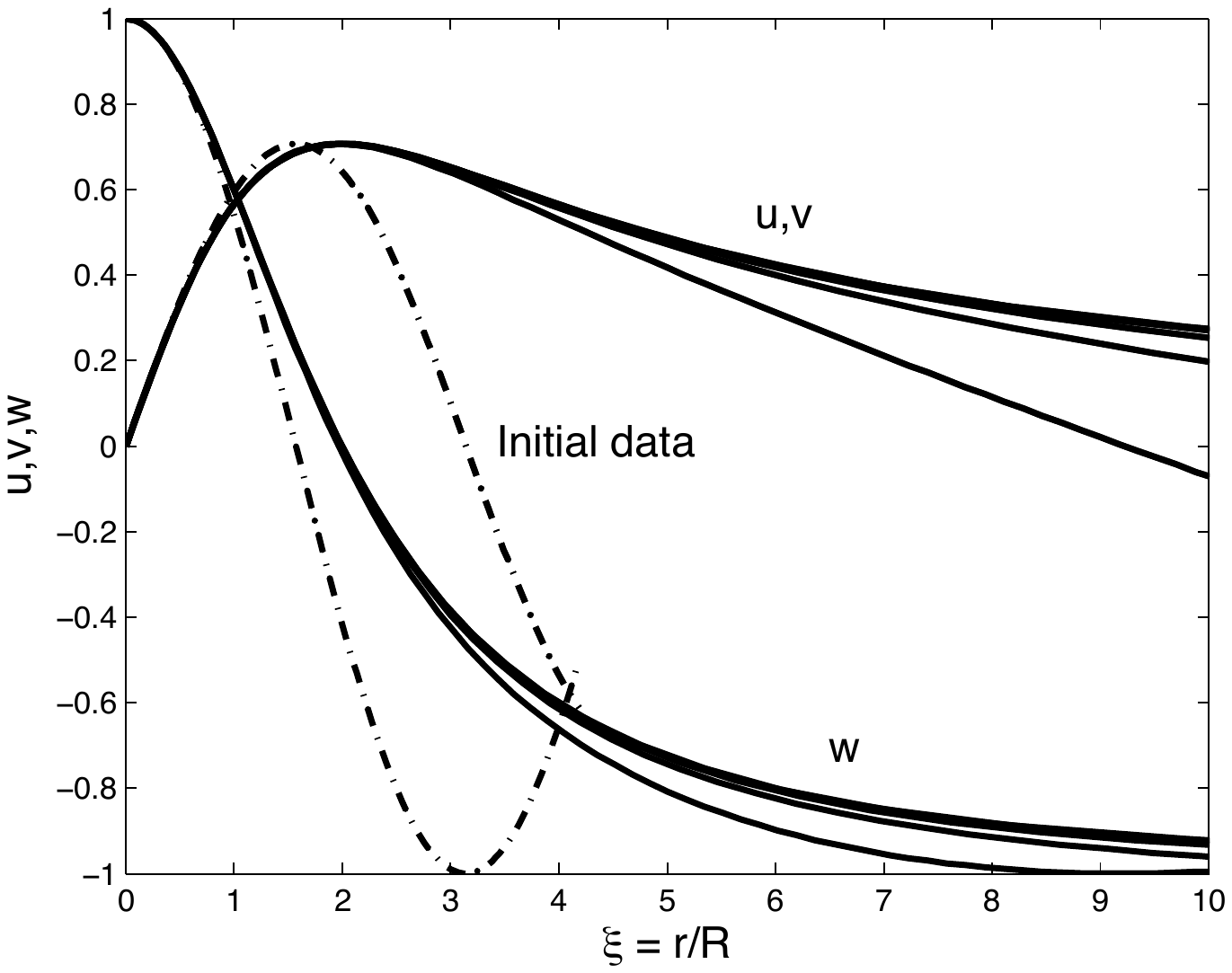}}}
\caption{Blowup of initial data (\ref{ThetaIC}) computed using (\ref{LLG_3comp}). Left: $\theta = \arctan(\sqrt{u^2+v^2}/w)$. Note that it again converges to the rescaled arctangent profile.  The initial
data is not monotone in~$r$ as now there is also a rotation in $\phi$.
Right: Solutions $u,v,w$ over time.}
\label{fig:Example2a1}
\end{figure}

First we reconsider the example above but using equation (\ref{LLG_3comp}) with the harmonic map case $\alpha = 0$, $\beta = 1$. We take the same initial as (\ref{ThetaIC}) and set
\[ u(0,t) =  v(0,t) =  \sin(\tfrac{4}{3} \pi r) / \sqrt{2}, \quad \mbox{and} \quad w(0,t) = \cos(\tfrac{4}{3} \pi r).\]
In Figures \ref{fig:Example2a} and \ref{fig:Example2a1} we see the same behaviour as observed above.  This is not surprising 
but a reassuring test of the numerics.

\begin{figure}[!t]
\centerline{{\includegraphics[width=.5\textwidth]{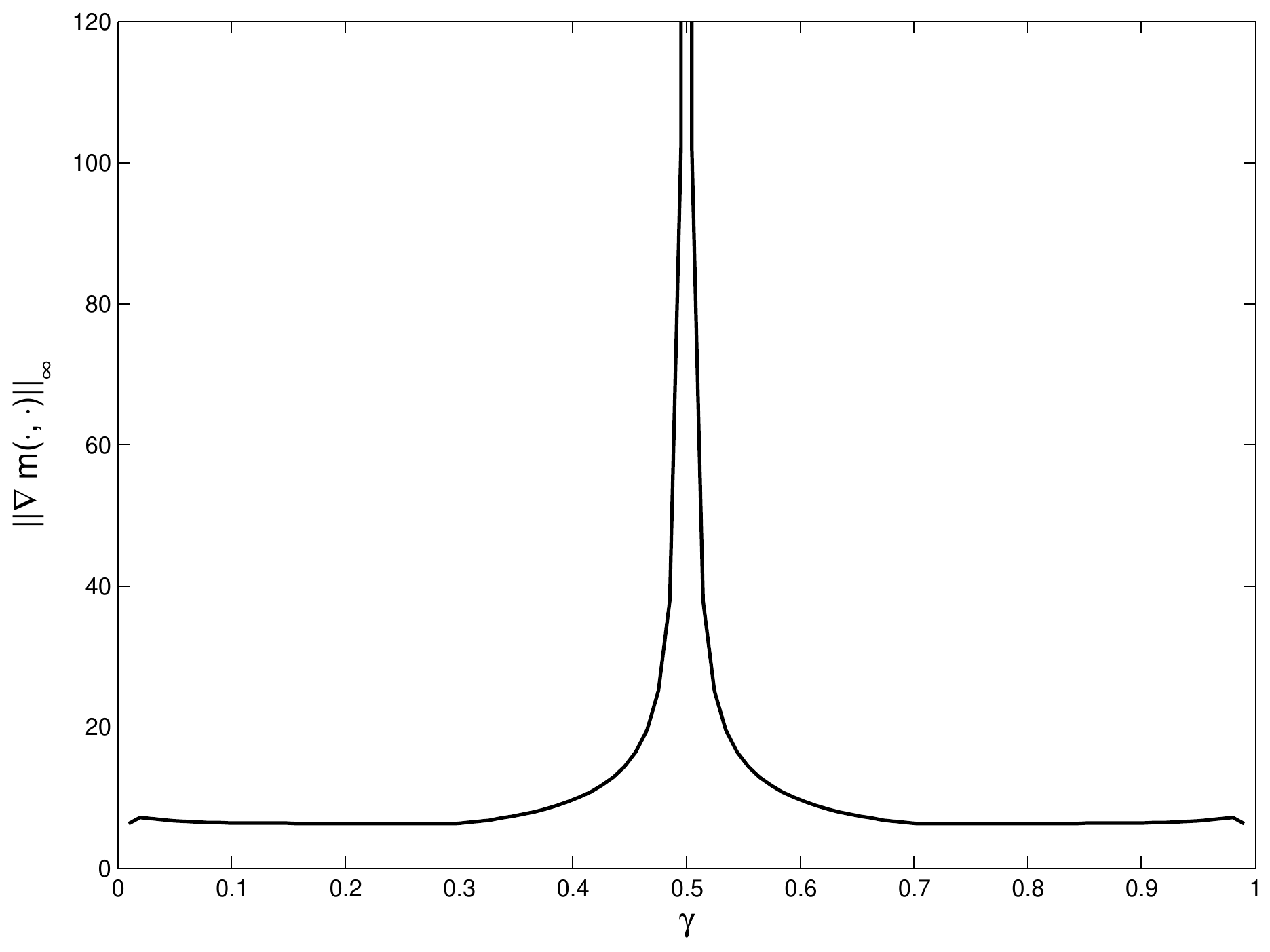}}
{\includegraphics[width=.495\textwidth]{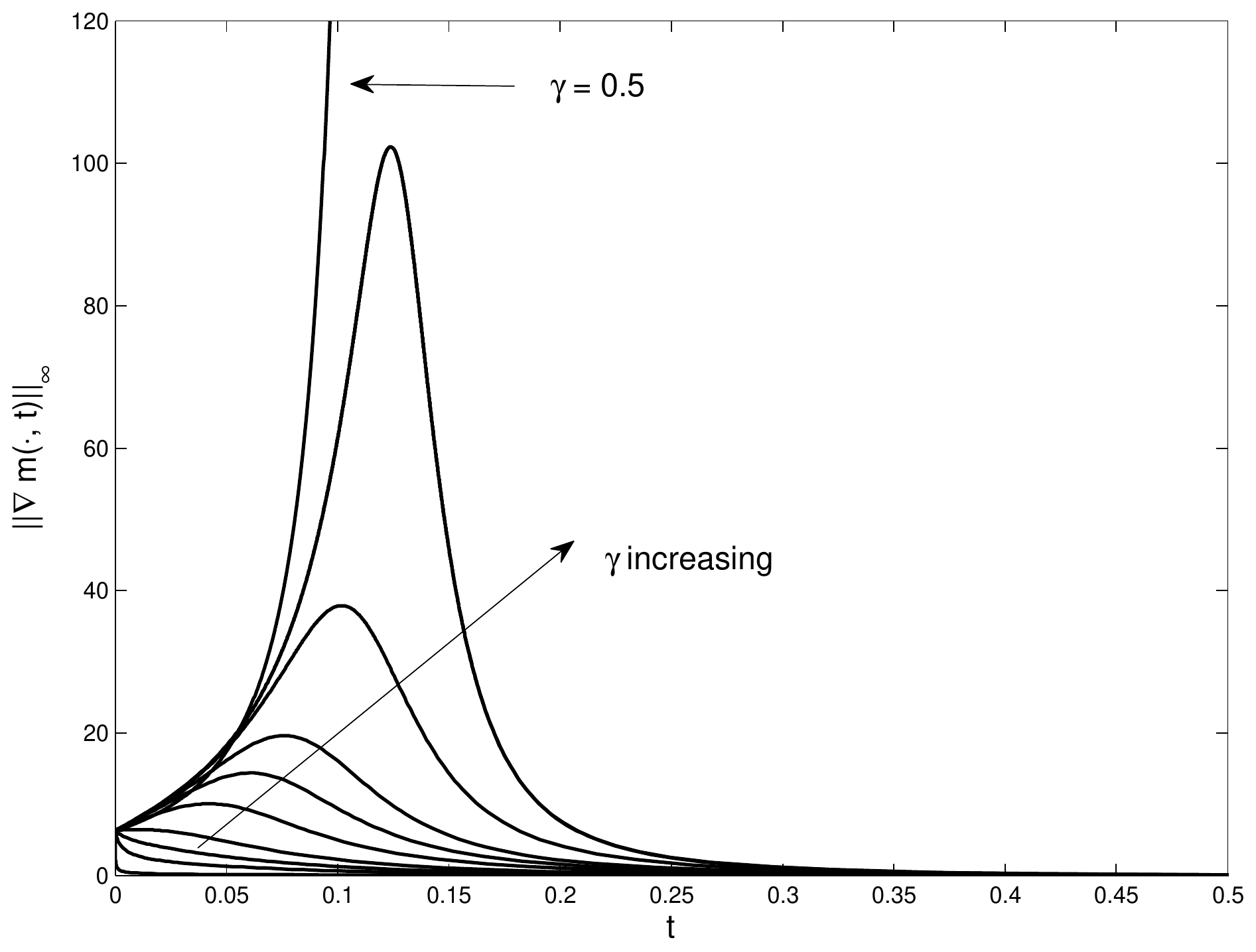}}}
\caption{(Left) $\|\nabla m\|_\infty$ as a function of $\gamma$ (for $\gamma = 0.5$ the computation was 
stopped when $\|\nabla m(\cdot, t)\|_\infty = 1e8$.) (Right) Growth and decay of $\|\nabla m(\cdot, t)\|_\infty$ over time for a sequence of values of $0 \le \gamma \le 0.5$ (in this case the dynamics are symmetric about $\gamma = 1/2$). }
\label{fig:Example2b}
\end{figure}

We now consider equation (\ref{LLG_3comp}) with $\alpha = 0$ and $\beta = 1$ for a family of initial data determined
via stereographic projection
\begin{equation}
 (u_\gamma,v_\gamma,w_\gamma) = \left(
   \frac{2x}{1+x^2+y^2},\frac{2y}{1+x^2+y^2},\frac{-1+x^2+y^2}{1+x^2+y^2} \right)
\label{GammaIC}
\end{equation}
where \[ x = \tan(-\pi/2 + r\pi), \quad\mbox{ and }\quad y = \tan(-\pi/2+\gamma \pi), \quad \mbox{for } \gamma \in [0,1], \]
which covers the sphere as $\gamma$ varies. From the discussion of Section \ref{sec:Top} we would expect blowup for a single value of $\gamma$ and decay to the 
stationary solution in all other cases. Figure \ref{fig:Example2b} shows $\max_{(r,t)} | \nabla m|$ as a function of the parameter $\gamma$ for a sequence of values of $\gamma$ and initial data (\ref{GammaIC}).

\begin{figure}[!tbh]
\centerline{{\includegraphics[width=.495\textwidth]{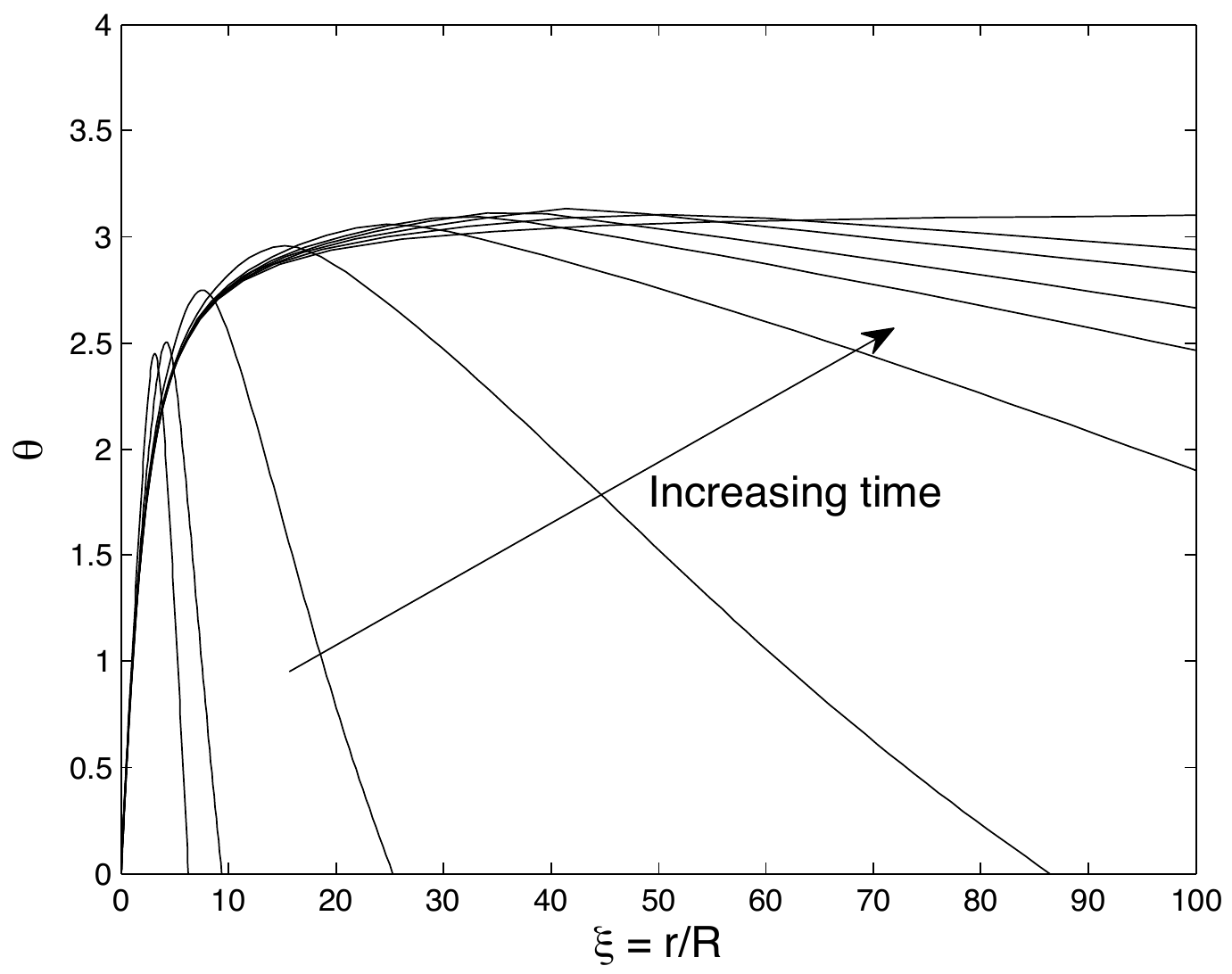}}
{\includegraphics[width=.5\textwidth]{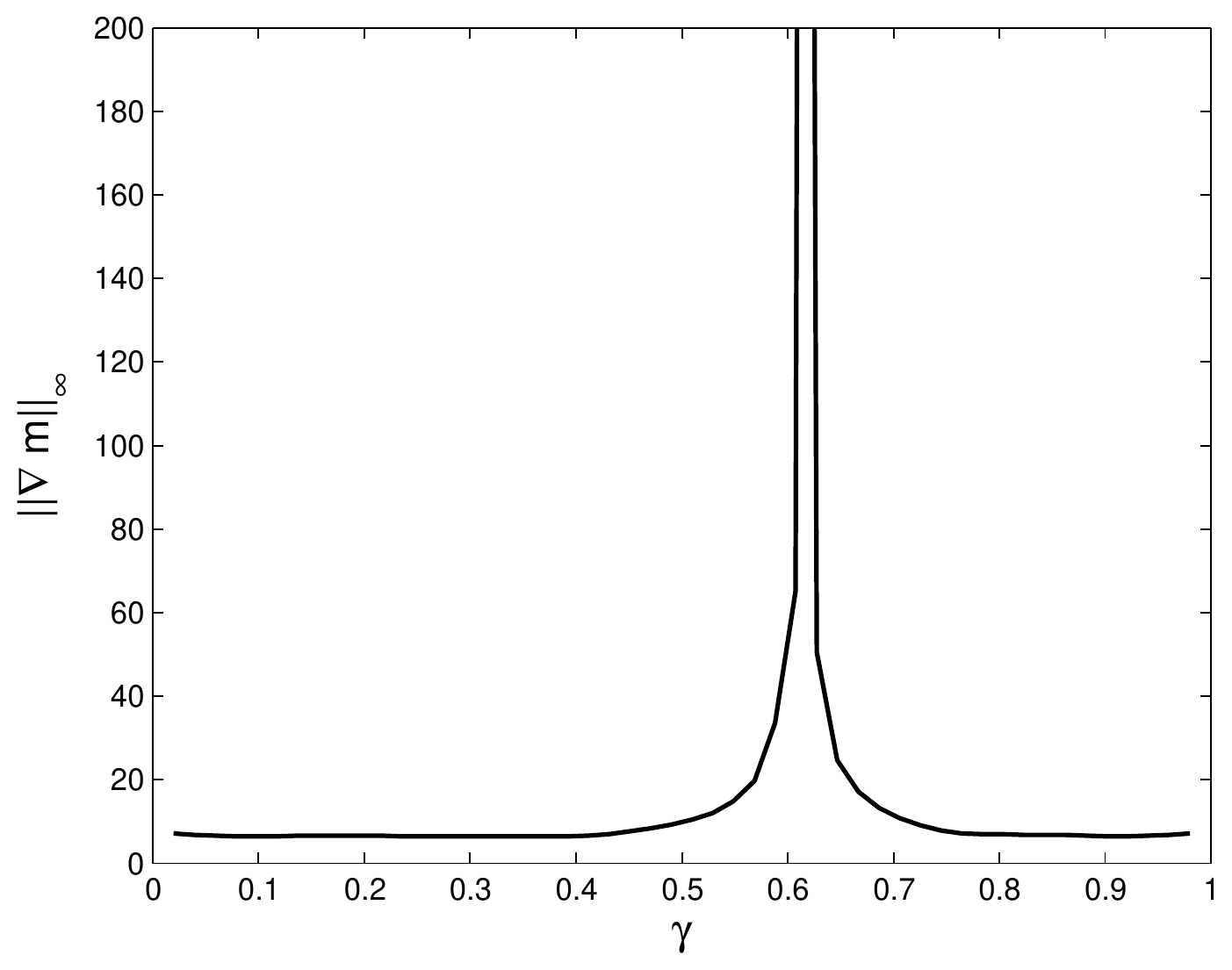}}}
\caption{(Left) Evolution of initial data (\ref{GammaIC}) with $\gamma = 0.612\ldots, \alpha = 1/\sqrt{2}$ and
$\beta = 1/\sqrt{2}$ (Right) Evolution of the maximum gradient for a sequence of values of $\gamma$ with $\alpha=\beta=1/\sqrt{2}$.}
\label{fig:Example3}
\end{figure}

{\bf Example 3 - Full Landau-Lifshitz-Gilbert $(\alpha > 0 )$}
We now consider the full Landau-Lifshitz-Gilbert equation with $\alpha \ge 0$ and
$\beta  = \sqrt{1-\alpha^2} \ge 0$. Figure \ref{fig:Example3} shows snapshots in time for $\alpha =1/\sqrt{2}$ and $ \beta =1/\sqrt{2} $ as well as $\max_{(r)} |\nabla m|$ over time for a sequence of values of $\gamma$ in ~(\ref{GammaIC}).
There is no qualitative difference to the case $\beta = 1, \alpha = 0$.
 
%
%

\begin{figure}[!tbh]
{\includegraphics[width=.32\textwidth]{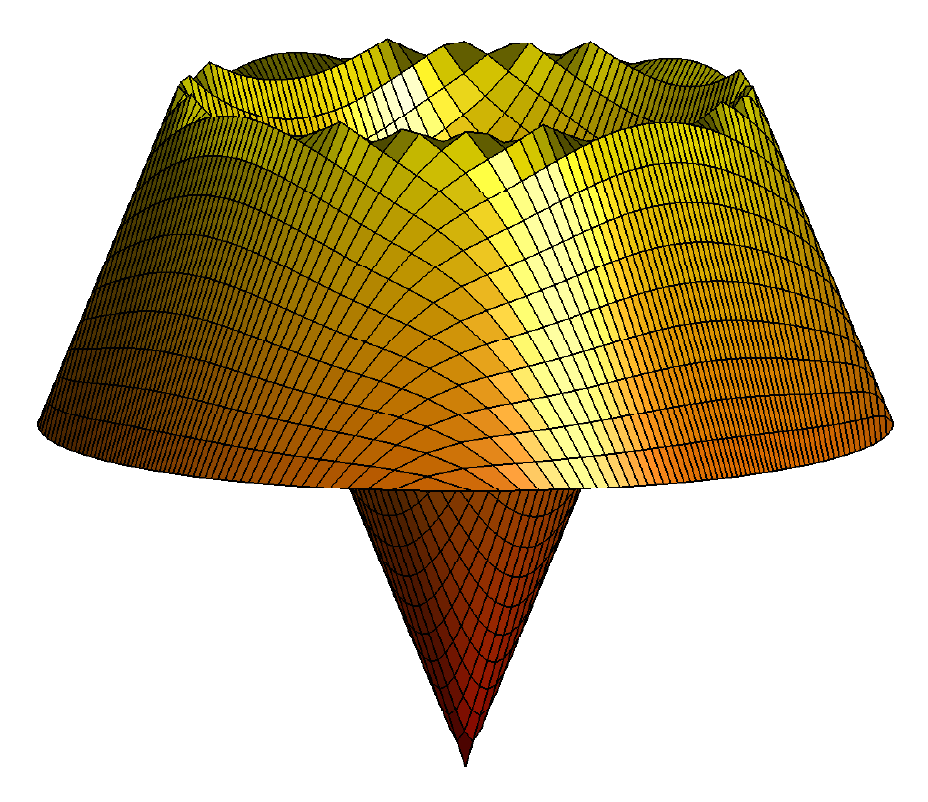}}
{\includegraphics[width=.32\textwidth]{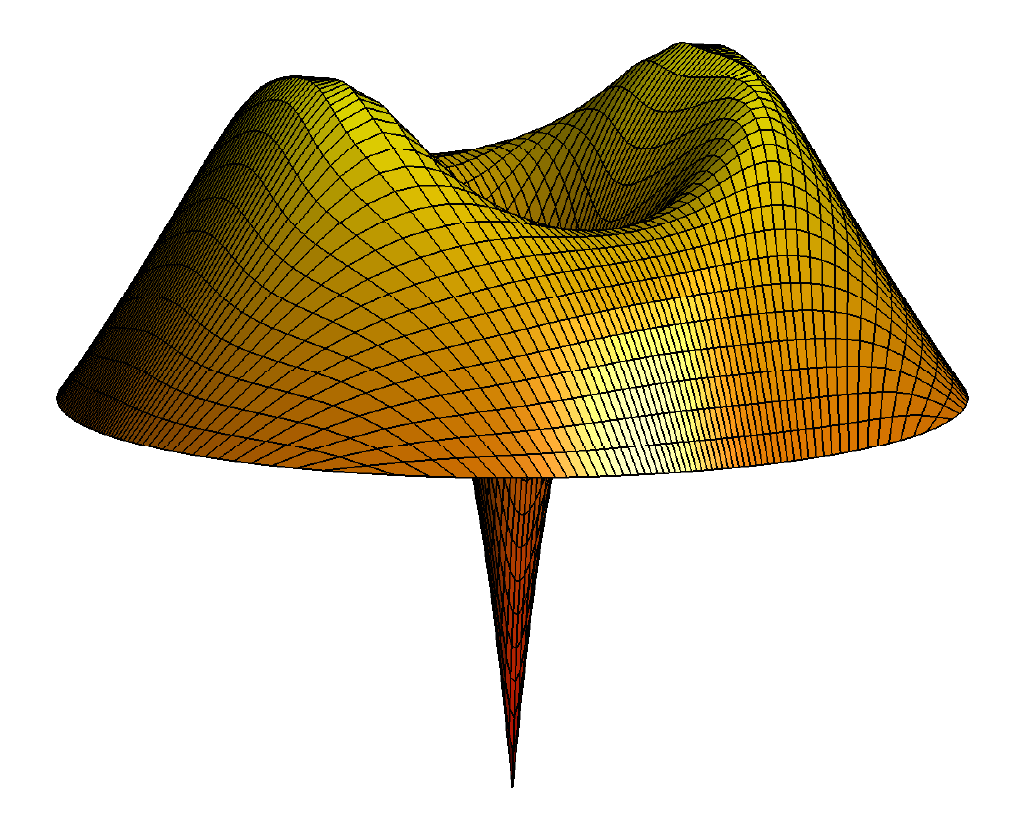}}
{\includegraphics[width=.32\textwidth]{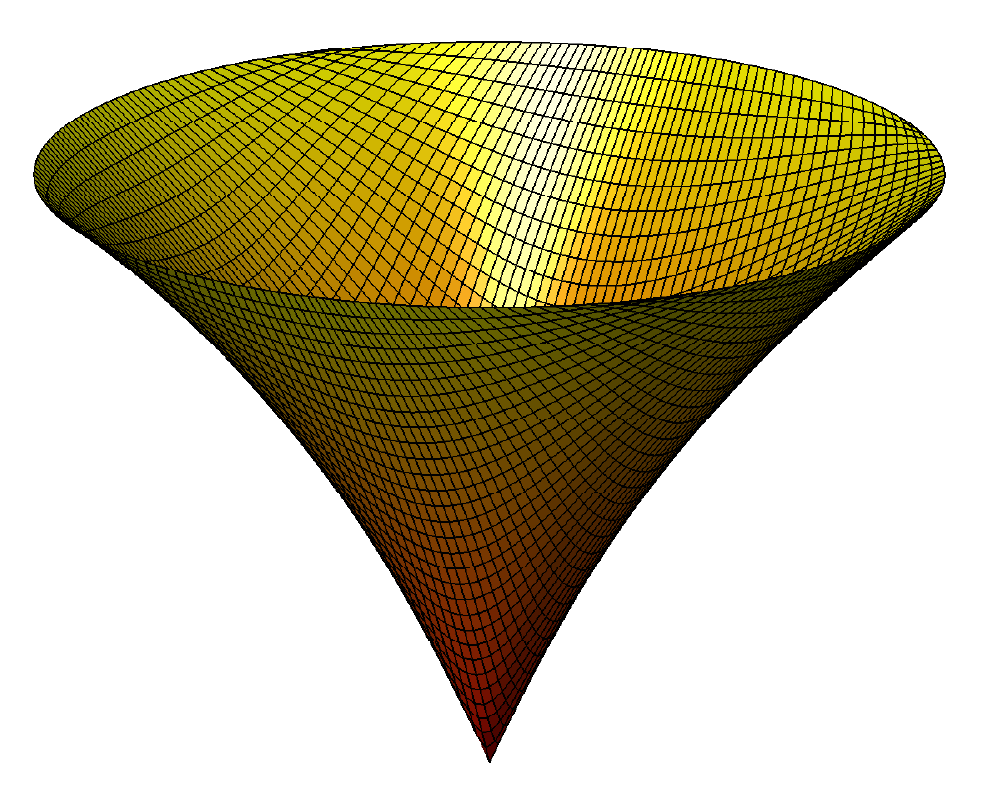}}
\caption{Evolution of the first component $m_1$ from non-radial initial data. The (Left) Initial data, $\| \nabla m \|_\infty = 23$. (Center) $\| \nabla m \|_\infty = 387$
(Right) $\| \nabla m \|_\infty = 12$.  Over time the asymmetry grows before the solution converges towards the
radially symmetric $\arctan$ profile.}
\label{fig:Example4a}
\end{figure}

\begin{figure}[!tbh]
{\includegraphics[width=.32\textwidth]{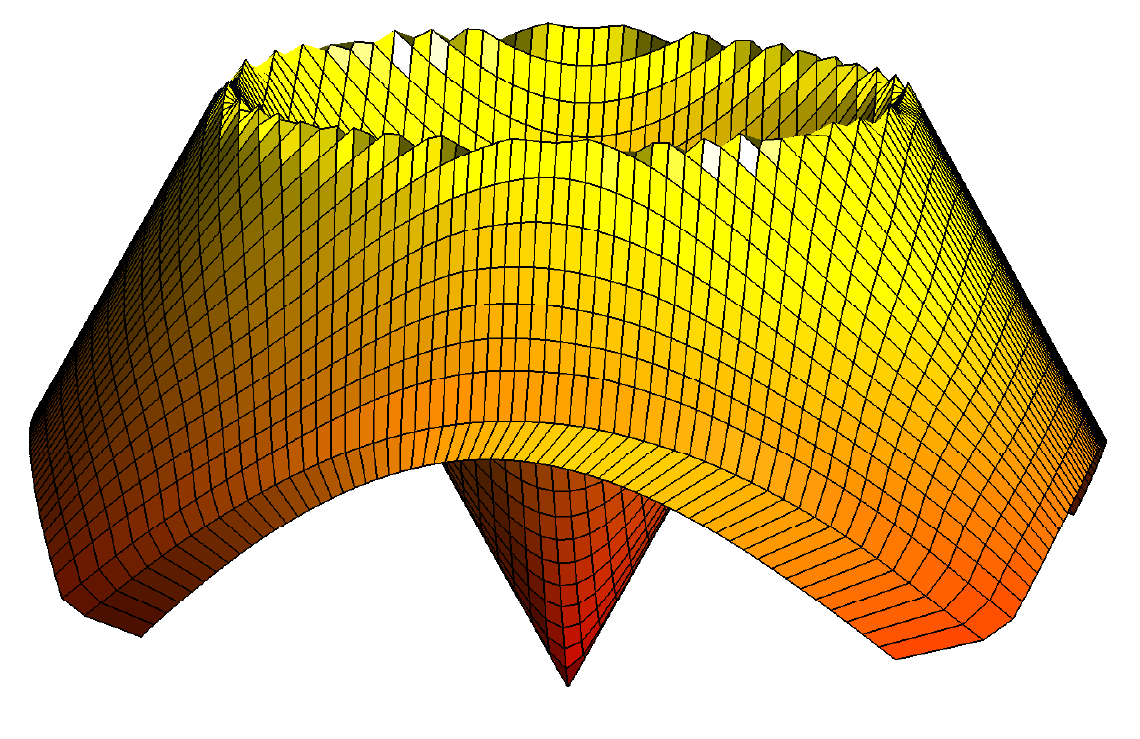}}
{\includegraphics[width=.32\textwidth]{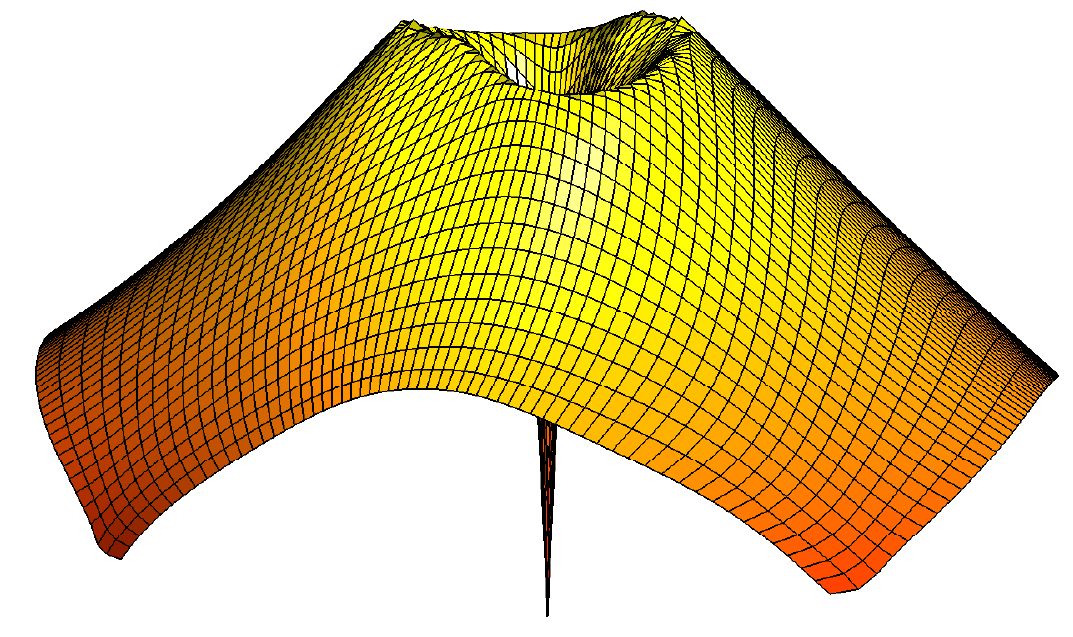}}
{\includegraphics[width=.32\textwidth]{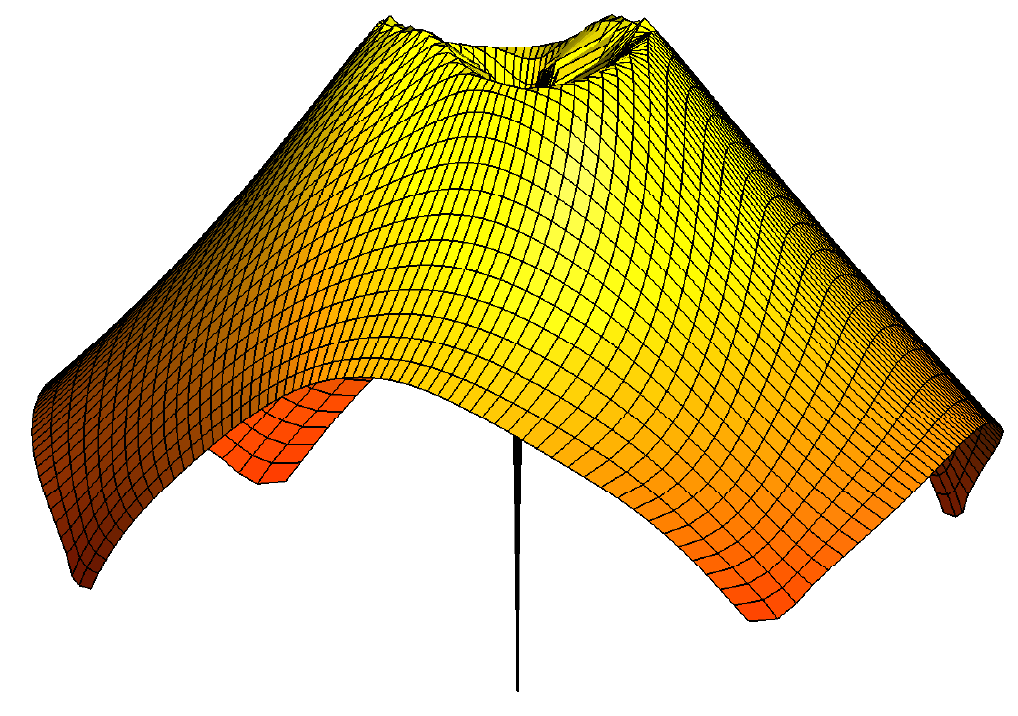}}
\caption{Evolution of the first component $m_1$ from non-radial initial data. (Left) Initial data, $\| \nabla m \|_\infty = 28$. (Center) $\| \nabla m \|_\infty = 953$
(Right) $\| \nabla m \|_\infty = 9.3e4$. }
\label{fig:Example4b}
\end{figure}

{\bf Example 4 - Full Landau-Lifshitz-Gilbert $(\alpha > 0 )$ in 2 dimensions}
We now consider the full Landau-Lifshitz-Gilbert equation with $\alpha > 0$ and
$\beta > 0$. Figure \ref{fig:Example4a} shows snapshots in time for $\alpha = 1/2, \beta = \sqrt{3}/2$ with initial data (\ref{GammaIC}) and $\gamma = 0.25$ and a small non-radial perturbation.
Figure \ref{fig:Example4b} shows snapshots in time for $\alpha = 1/2, \beta = \sqrt{3}/2$ but now we have taken
a larger non-radial perturbation of (\ref{GammaIC}) and varied $\gamma$ until we had evidence of blowup.  Here $\| \nabla  m \|_\infty = \max_{j=1\dots3} ((\partial_1 m_j)^2 + (\partial_2 m_j)^2)^{1/2}$
changes  almost four orders of magnitude before the computation halts.
   
Even though the analysis above is for radial initial data we can find solutions that lead to blowup with carefully tuned parameters
specific to given non-radial initial data. This is not necessarily a true blowup solution but rather a \emph{numerical} one in the sense that it focusses to such a degree that we cannot continue the computation.

\section{Behavior of near-blowup solutions}
\label{sec:MoreAsymp}

Instability leads to a reconfiguration described by a quick rotation of the
sphere that had almost bubbled off. The derivation and asymptotics of this quick rotation are presented in Section~\ref{sec:rotation} below. This is highly relevant for the problem of continuing the exceptional solution that \emph{does} blowup after its blowup time, as explained in Section~\ref{sec:continuation}. 

\subsection{The quick rotation}
\label{sec:rotation}

Here we present the asymptotics of near-blowup solutions.
The inner scale (in the domain; it describes a sphere in
the image, or a semicircle in equivariant coordinates) 
is given by the usual $\xi = r/R(t)$ with
\[
  \theta \sim 2 \arctan \xi + (\beta R'R- \alpha R^2 C')[\xi-\xi\ln \xi]
  \qquad \text{for large } \xi,
\]
and
\[
  \phi = C(t) + (\beta R^2 C'+ \alpha  R'R) [ \frac{1}{2} \xi^2 (\ln \xi -1)]
  \qquad \text{for large } \xi.
\]
For the outer scale (representing a small neighborhood of the south pole $S$ in the image)
we introduce a fast time scale $t=T+\eps^2 \tilde{t}$. On this time scale the dynamics takes place at a small spatial
scale $r=\eps\tilde{r}$, 
but large compared to $R(t)$, i.e., $R \ll \eps \ll 1$,
where the solution is described by ($z$ representing coordinates in the
tangent plane at the south pole as in Section~\ref{sec:outer})
\[
  z_{\ttt} = (\beta-\alpha i) \left( z_{\tr\tr} +\frac{1}{\tr}z_{\tr}  -
  \frac{1}{\tr^2}z \right),
\]
with solution ($\sigma= \sigma_r + i \sigma_i$ and
$\gamma=\gamma_r+i\gamma_i$)  
\[
  z = \sigma(\ttt) \tr^{-1} + \gamma(\ttt) \tr + \dots .
\]
Looking at the modulus and argument of $z$ we obtain for small $r$
\[
  |z| \sim \sqrt{\sigma_r^2+\sigma_i^2} r^{-1} +
  \frac{\sigma_r\gamma_r+\sigma_i\gamma_i}{\sqrt{\sigma_r^2+\sigma_i^2}} r ,
\] 
and
\[
  \text{arg} z = \arctan \frac{\sigma_i}{\sigma_r} +
  \frac{\sigma_r\gamma_i-\sigma_i\gamma_r}{\sigma_i^2+\sigma_r^2} r^{2}.
\]
Matching $|z|$ to $\pi -\theta$ and $\arg z$ to $\phi$, the matching conditions read
\begin{alignat*}{3}
|z| &:\quad &\tr^{-1}&:~&   2 \eps^{-1} R &\sim (\sigma_r^2+\sigma_i^2)^{1/2} , \\
   & &\tr^{1}&:& -\eps^{-1} (\beta R' - \alpha R C') \ln R  
&\sim \frac{\sigma_r\gamma_r+\sigma_i\gamma_i}{(\sigma_r^2+\sigma_i^2)^{1/2}} , \\
\text{arg}z &: &\tr^{0}&:& C &\sim \arctan \frac{\sigma_i}{\sigma_r} ,\\
& &\tr^{2}&:& - \frac{1}{2}(\beta C' +\alpha R' R^{-1}) \ln R & \sim 
\frac{\sigma_r\gamma_i-\sigma_i\gamma_r}{\sigma_r^2+\sigma_i^2} .
\end{alignat*}
In the remote region we have $z(r) \sim q r$ for small $r$ for 
some $q\in \mathbb{C}$, where $q \sim \sigma \sim k_2/\tau$ is small close to blowup, as explained in
Section~\ref{sec:Asymp}. 
And as before, by rotating the sphere we may assume that $q=-q_0$, with $q_0>0$ real.
By matching it follows that $z(\tr) \sim - q_0 \eps \tr $ for large $\tr$,
hence $\gamma_r \approx - \eps q_0$ and $\gamma_i \approx 0$.

Hence, by rearranging the terms we obtain 
\begin{alignat*}{1}
  2 \eps^{-1} R &\sim (\sigma_r^2+\sigma_i^2)^{1/2} , \\
  C &\sim \arctan \frac{\sigma_i}{\sigma_r} , \\
  \eps^{-1}  R' \ln (1/R)   
&\sim \frac{-\beta \sigma_r q_0\eps +\alpha\sigma_i q_0\eps}{(\sigma_r^2+\sigma_i^2)^{1/2}} , \\
 \eps^{-1}  C' R \ln (1/R)  
&\sim \frac{\beta \sigma_i q_0\eps + \alpha \sigma_r q_0\eps}{(\sigma_r^2+\sigma_i^2)^{1/2}} .
\end{alignat*}
Let us again remove $\alpha$ and $\beta$ from the formulas by 
setting $\mu_r = \alpha \sigma_r - \beta \sigma_i$ and $\mu_i=\beta\sigma_r +\alpha
\sigma_i$. In complex notation: $\mu_r+\mu_i i = (\beta+ \alpha
i)(\sigma_r+\sigma_i i)$. 
Furthermore, write $\tC=C+\arctan \frac{\alpha}{\beta}$. This leads to
\begin{alignat*}{1}
  2 \eps^{-1} R &\sim (\mu_r^2+\mu_i^2)^{1/2} , \\
  \tC &\sim \arctan \frac{\mu_i}{\mu_r} , \\
  \eps^{-2}  R' \ln(1/R)  
&\sim -\frac{\mu_r q_0 }{(\mu_r^2+\mu_i^2)^{1/2}} ,\\
 \eps^{-2} R\tC' \ln(1/R) 
&\sim \frac{\mu_i q_0 }{(\mu_r^2+\mu_i^2)^{1/2}}  .
\end{alignat*}

\begin{figure}
\centerline{\includegraphics{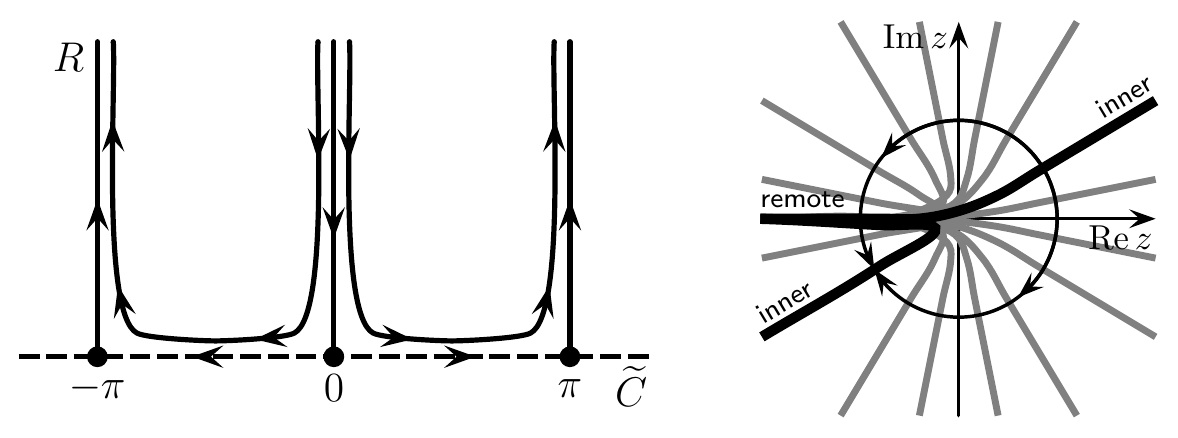}}
\caption{Left: in the $(\tC,R)$ phase plane the ``stable manifold'' of the blowup point acts as separatrix. Right: geometrically it is the boundary between a rotation over an angle $\pi$ or $-\pi$.}
\label{fig:separatrix}
\end{figure}

Looking at the matching conditions, 
we write $\mu_r = 2\eps^{-1} R \cos \tC$ and $\mu_i = 2\eps^{-1} R
\sin \tC$, with $\frac{dR}{d\ttt} = O(\eps^2)$, which we can neglect on this time scale.
We are left with the 
dynamics of $\tC$, determined by the remaining equation
\[
  \frac{d\tC}{d\ttt} = \frac{q_0 \eps^2}{R \ln (1/R)}  \sin \tC.
\]
We see that the correct time scale is $\eps^2 = \frac{R \ln (1/R)}{q_0}$,
which is smaller the closer we are to blowup (and the larger $q_0$ is). Notice  that indeed $\eps \gg R$ since $q_0 = O(1/\ln R)$ near blowup as discussed before, demonstrating self-consistent separation of spatial scales.  
The angle $\tC$ thus approaches $\pm \pi$ depending on the initial data,
unless $\tC=0$. The quick rotation due to the instability is described
by $\frac{d\tC}{d\ttt} = \sin \tC$, and the solution, in the original time variable, is $\tC(t)=
\pm [\frac{\pi}{2}+\arctan[\sinh(\eps^{-2}(t-T)+c_0))]$, with $c_0\in\mathbb{R}$. This shows that the blowup solution acts as a separatrix between rotations in two opposite directions, see Figure~\ref{fig:separatrix}.
The dependence on $\alpha$ and $\beta$ in this scale is only through the fixed rotation $C=\tC-\arctan \frac{\alpha}{\beta}$. 

\subsection{Numerical investigation of near blowup}

\begin{figure}[!t]
{\includegraphics[width=.32\textwidth]{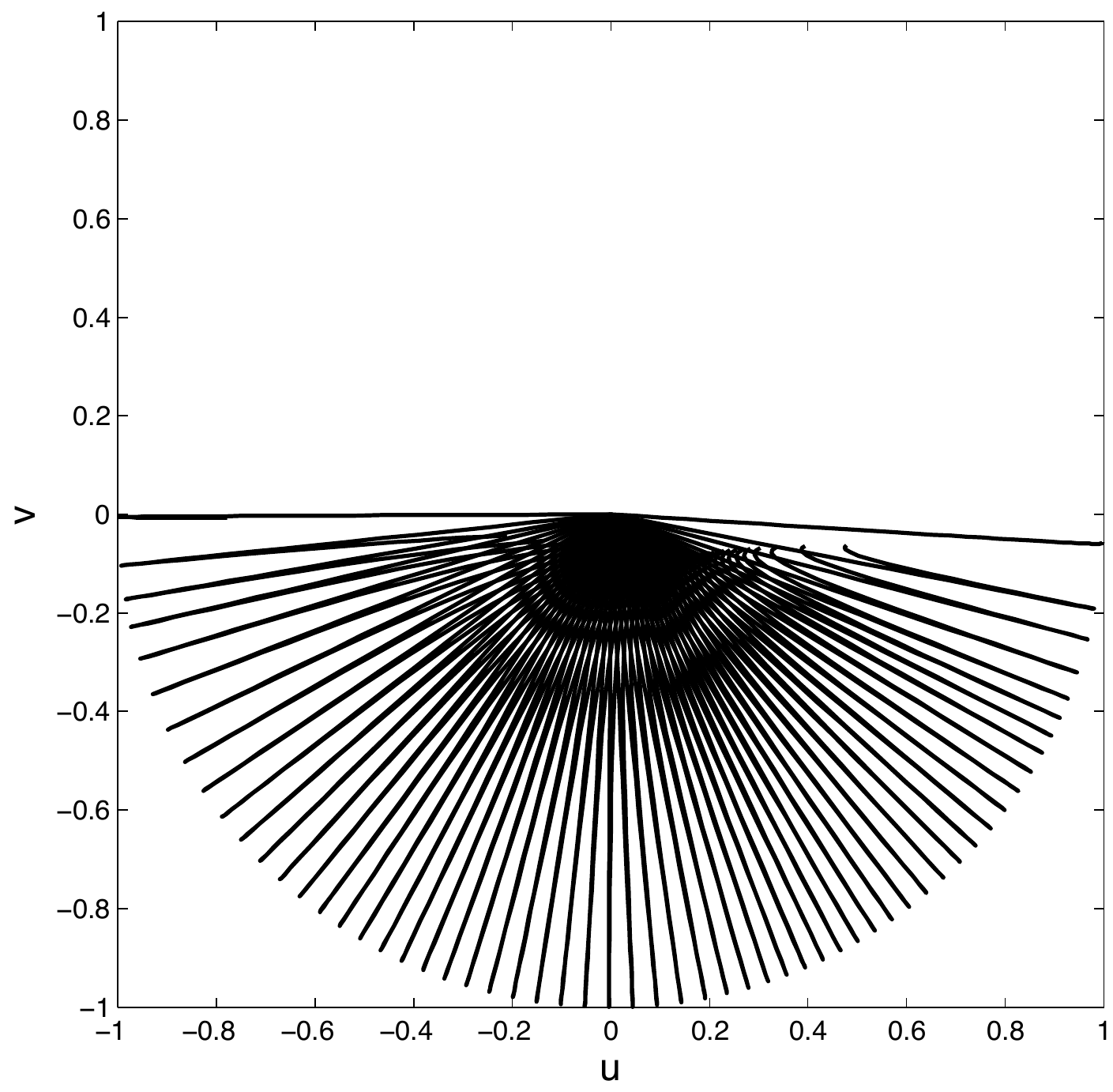}}
{\includegraphics[width=.32\textwidth]{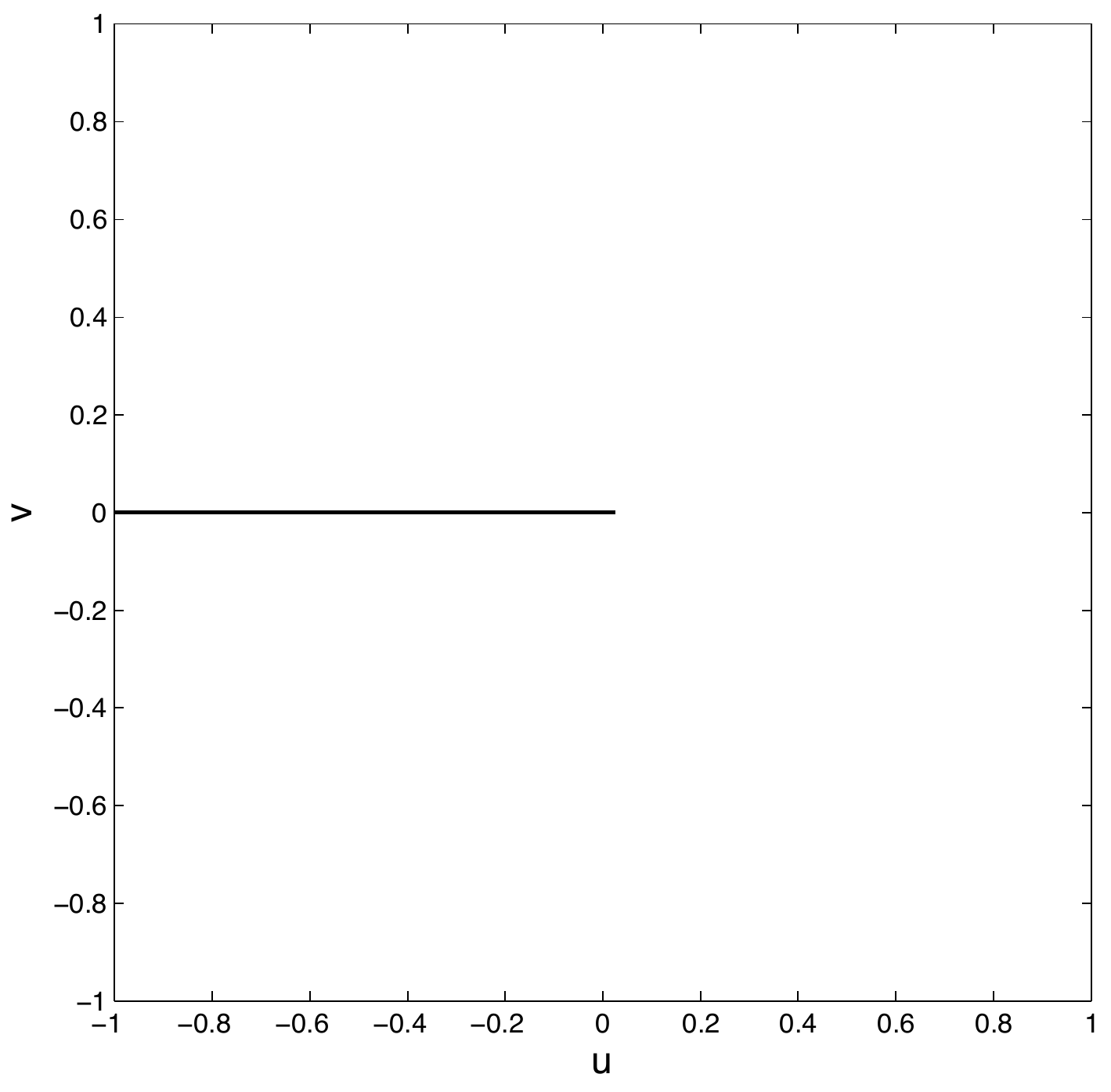}}
{\includegraphics[width=.32\textwidth]{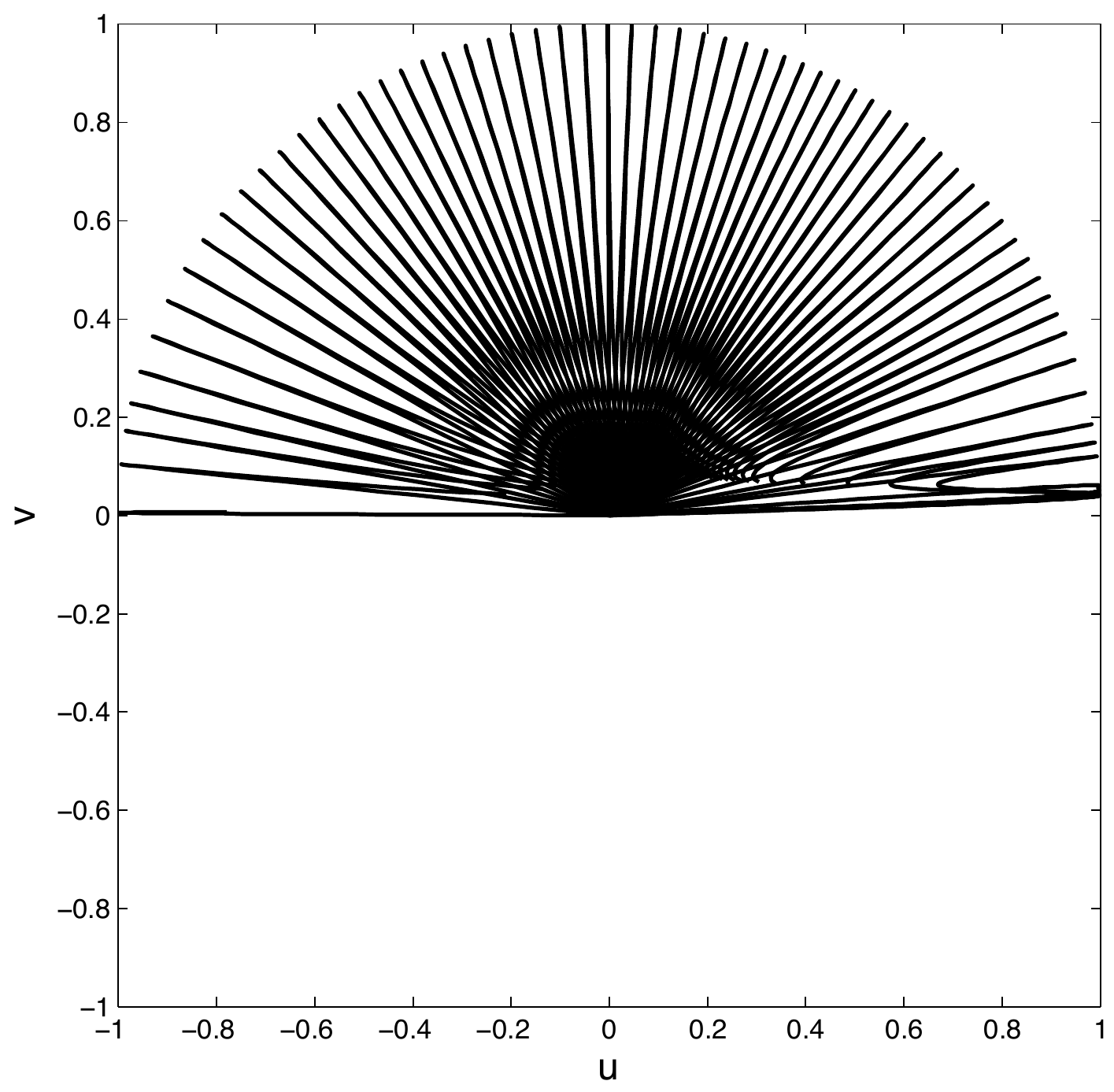}}
\caption{(Left) $\gamma < .5$. (Centre) $\gamma = 0.5$. (Right) $\gamma > .5$.  In this sequence we plainly see the role of the blowup solution as a separatrix.  There is blowup for $\gamma = 0.5$ and rotation away form blowup in opposite
directions for $\gamma < 0.5$ than for $\gamma > 0.5$. Here $u$ and $v$ have been plotted in the spatiotemporal regime 
close to blowup.}
\label{fig:Example5}
\end{figure}

In the previous Section we saw solutions whose gradient grew dramatically and then decayed as well as those that show blowup.  We can investigate the near blowup solutions in the context of the previous subsection by plotting $\phi = \arctan(v/u)$ in the region where the norm is large and also by plotting the dynamics in the $(u,v)$-plane. Figure 
\ref{fig:Example5} shows the latter for three runs with $\alpha =0$ and $\beta = 1$ for three values of $\gamma$ near
the critical value $\gamma = 0.5$.  In the two cases with $\gamma \ne 1/2$ we see the initial motion towards
the singularity followed by decay to a regular equilibrium whereas $\gamma = 1/2$ leads to the separatrix blowup behaviour.

\subsection{Continuation after blowup}
\label{sec:continuation}

Starting from smooth initial data, solutions to (\ref{LLG1}) are unique as long as they are classical, and finite time blowup may indeed occur for the (radially symmetric) harmonic map heatflow~\cite{CDY}. At a blowup point (in time) the strong solution terminates (at least temporarily). 
On the other hand, it is known that weak solutions of (\ref{LLG1}) exist globally in time~\cite{Struwe2,Struwe,AlougesSoyeur,GuoHong,BertschPodioValente}.
It is well established that such weak solutions are not unique~\cite{AlougesSoyeur,Coron,BertschDalPassoHout,Topping}.
It is thus of interest to come up with criteria that select the ``most appropriate'' weak solution. In other words, how should one continue a solution after the blowup time?

For the harmonic map heatflow in two dimensions it has been shown~\cite{Struwe2,Freire} that one uniqueness criterion is non-increasing energy
\[
  e(t) = \frac{1}{2} \int_{D^2} |\nabla m(t)|^2 , 
\]
i.e., there is exactly one weak solution that has non-increasing energy $e(t)$ for all $t \in [0,\infty)$.  
Furthermore, the energy of a solution jumps down by at least $4\pi$ at a singularity.

The co-dimension one character of blowup and our analysis of near-blowup solutions in Section~\ref{sec:rotation} leads us to propose a different scenario for continuation after blowup. It has the important advantage of \emph{continuous dependence on initial data} for times after blowup.
The scenario is identical for all parameter values $\alpha$ and $\beta$.
Namely, consider an equivariant solution of (\ref{LLG1}) that blows up as $t \uparrow T$. The blowup behaviour is characterized by a length scale $R(t)\to 0$ as $t \uparrow T$ and
\[
  \phi(r,t) \to \overline{\phi}   \quad\text{and}\quad
  \theta(r,t) \sim 2 \arctan \frac{r}{R(t)}
  \qquad\text{for } r=O(R(t)) \text{ as } t \uparrow T,
\]
i.e., geometrically speaking a sphere bubbles off at $t=T$.
Based on the analysis of near-blowup solutions, we propose to continue the solution for $t>T$ by \emph{immediately} re-attaching the sphere, \emph{rotated over an angle $\pi$} with respect to the bubbled-off sphere:
\[
  \phi(r,t) \to \overline{\phi}+ \pi  \qquad \text{as } t \downarrow T,
\]
with $\theta(r,t) \sim 2 \arctan \frac{r}{\widetilde{R}(t)}$ for $r=O(\widetilde{R}(t))$, and $\widetilde{R}(t) \to 0$ as $t \downarrow T$.
By rotating the re-attached sphere (also referred to as a reverse bubble \cite{Topping}) over an angle $\pi$, this continuation framework leads to continuous dependence on initial data, since nearby solution that avoid blowup undergo a rapid  rotation over an angle $\pi$, as derived in Section~\ref{sec:rotation}. 

A solution that is continued past blowup through the re-attachment of a rotated sphere, does not have a monotonically decreasing energy $e(t)$. However, the renormalized energy 
\[
  \overline{e}(t) = \left\{ \begin{array}{ll}
  e(t)& \text{for }  t \neq T, \\
  e(T)+4\pi \quad & \text{for }  t=T,
  \end{array} \right.
\]
is continuous and decreases monotonically.

In the radially symmetric harmonic map heatflow case described by (\ref{radHM})
this scenario corresponds to 
\[
 \theta(0,t)=  \left\{ \begin{array}{ll}  
  0 \quad& \text{for } t <T,\\
  \pi& \text{for } t =T,\\
  2\pi & \text{for } t >T.
\end{array} \right.
\]
For this particular case it has been proved~\cite{vanderHout} that such a re-attachment leads to a \emph{unique} solution for $t>T$. For general equivariant solutions of (\ref{LLG1}) such an assertion remains an open problem. 
Moreover, all conclusions in Sections~\ref{sec:Asymp} and~\ref{sec:MoreAsymp}, as they follows from formal matched asymptotics, require mathematically rigorous justification.

\section{Other $n \geq 2$}
\label{sec:AsympN}

We now summarize the calculations for $n=2,3,\dots$, 
following the same methodology as for $n=1$,  but now the formulas are simpler (since only the
first term in the expansion in the outer scale is needed).
As analyzed in~\cite{BHK}, for $n\geq 2 $ blowup is in infinite time.
We shall only consider blowup with one sphere bubbling off (i.e.\ $\theta_1
\in (\pi,2\pi]$) in the harmonic map flow case; the adaptation to the general case $(\alpha \neq 0$) is analogous to Sections~\ref{sec:Asymp} and~\ref{sec:MoreAsymp}.  
The $\theta$-component of the large $\xi$ asymptotics in the inner scale 
was already calculated in~\cite{BHK}:
\[
  \theta \sim \pi - 2  \xi^{-n} + (\beta R'R- \frac{\alpha}{n} R^2 C') 
  \left(\frac{n}{2n-2} \xi^{-n+2}-E_n \xi^n\right) ,
\] 
with 
\[
  E_n =\int_0^\infty \frac{s^{2n+1}}{(1+s^{2n})^2} ds
  = \frac{\pi}{2n^2 \sin \frac{\pi}{n}} \,.
\]
With  $n\geq 2$ equation~(\ref{e:phi1}) for $\phi_1$ is now replaced by
\[
  \phi_{1\xi\xi} +
  \frac{(2n+1)-(2n-1)\xi^{2n}}{\xi(1+\xi^{2n})}\phi_{1\xi}=1.
\]
Using the boundary condition $\phi_{1\xi}(0)=0$, we find
\[
  \phi_{1\xi} = \frac{\xi^{-2n}+2+\xi^{2n}}{\xi} 
  \int_0^\xi \frac{s^{2n+1}}{(1+s^{2n})^2} ds,
\]
which has asymptotic behaviour $\phi_{1\xi} \sim E_n \xi^{2n-1} -
\frac{1}{2n-2} \xi$ as $\xi \to\infty$. We thus find that
\[
  \phi \sim C + (\beta R^2 C'+ \alpha n R'R)
  \left(\frac{E_n}{2n} \xi^{2n} - \frac{1}{4(n-1)} \xi^2\right).
\]
Since the blowup for $n \geq 2$ occurs as $t \to \infty$,
the outer variables are just the $\cO(1)$ $t$ and $r$ (i.e. not self-similar), and the equation becomes 
\[
  z_t = (\beta-i\alpha) \left(z_{rr}+\frac{1}{r} z_r - \frac{1}{r^2} z \right).
\]
The solution is asymptotically given by (with $\gamma$ and $\sigma$ complex valued) 
\[
  z = \gamma(t) r^n + \sigma(t) r^{-n} + \dots
\]  
This solution needs to match into the remote region where $\phi=\pi$,
$\theta= \pi- 2\arctan q_0 r^n$, with $q_0= \tan \frac{\pi-\theta_b}{2}$. 
Hence $\gamma \approx -q_0 \in \mathbb{R}$. 

This leads to the matching conditions (see also Section~\ref{sec:rotation})
\begin{alignat*}{3}
|z| &:\quad &r^{-n}&:~&   2 R^n &\sim (\sigma_r^2+\sigma_i^2)^{1/2} , \\
   & &r^{n}&:&  (\beta R'R^{1-n} - \frac{\alpha}{n}R^{2-n}C')E_n  
   &\sim -\frac{\sigma_rq_0}{(\sigma_r^2+\sigma_i^2)^{1/2}} ,  \\
\text{arg}z &: &r^{0}&:& C &\sim \arctan \frac{\sigma_i}{\sigma_r} , \\
& &r^{2n}&:& \frac{1}{2n}(\beta R^{2-2n}C' +\alpha n R' R^{1-2n})E_n &\sim 
\frac{\sigma_iq_0}{\sigma_r^2+\sigma_i^2} .
\end{alignat*}
As before, let us transform the equation to remove the explicit dependence on $\alpha$ and $\beta$ by
setting $\lambda_r+\lambda_i i = (\beta+ \alpha
i)(\sigma_r+\sigma_i)$. 
Furthermore, write $\tC=C+\arctan \frac{\alpha}{\beta}$ to obtain
\begin{alignat*}{1}
  2R^n &\sim |\lambda|,\\
  \tC &\sim \text{arg} \lambda,\\
  2R^2\tC' E_n &\sim n \lambda_i q_0 ,\\
  2RR' E_n &\sim - \lambda_r q_0 .
\end{alignat*}
Hence $\lambda_r=2R^n \cos \tC$ and $\lambda_i = 2R^n \sin \tC$ and the
remaining system is
\begin{alignat*}{1}
  \tC' = \frac{nq_0}{E_n} R^{n-2}\sin \tC,\\
  R' = - \frac{q_0}{E_n} R^{n-1} \cos \tC, 
\end{alignat*}
from which we easily conclude that blowup is unstable,
since the equilibria $R=0$ and $\tC=k\pi$ are all unstable
(in the $\tC$ direction if $k$ is even, and in the $R$ direction if $k$ is odd).

\section{Conclusions}
\label{sec:Concl}

In this paper we have clearly demonstrated that blowup in the full Landau-Lifshitz-Gilbert equation is possible but that it is not generic.  Instead, we have identified finite-time blowup as a co-dimension one phenomenon possible only for specially chosen initial data. It is analogous to a saddle point along whose unstable manifold the flow is much slower than on the stable one. This means that while actual finite time blowup occurs for initial data on a set of measure zero, there is a wide set of initial data for which the solution gradient does increase significantly and may appear to blow up in numerical simulation.  

While we agree with the results in \cite{BKP, KP} about {\em discrete}  blowup in this equation, their computations do not indicate generic blowup in the continuous problem.  Blowup in this problem corresponds to energy concentration at small scales and so will vanish on any fixed grid with limited resolution. The numerical results in \cite{BKP} show changes in energy of little more than one order of magnitude and are resolution limited with $h=1/64$. Instead of continuous blowup, growth in those results halt when the solution can no longer be resolved and the authors carefully chose to call that {\em discrete} blowup \cite{BKP}.

In many other problems this would not be an issue as blowup typically occurs only because some small scale physical effects (surface-tension, high-order diffusion, saturation, etc.) have been neglected.  In those cases blowup means loss of model validity.  However, in this problem, it is the geometry of the target manifold that leads to the singularity and {\em it cannot be avoided by simply adding a regularizing term}.

\bibliographystyle{alpha}
\bibliography{refs}

\begin{thebibliography}{10}

\bibitem{AlougesSoyeur}
{\sc F.~Alouges and A.~Soyeur}, {\em On global weak solutions for
  {L}andau-{L}ifshitz equations: existence and nonuniqueness}, Nonlinear Anal.,
  18 (1992), pp.~1071--1084.

\bibitem{AngenentHulshof}
{\sc S.~Angenent and J.~Hulshof}, {\em Singularities at {$t=\infty$} in
  equivariant harmonic map flow}, in Geometric evolution equations, vol.~367 of
  Contemp. Math., Amer. Math. Soc., Providence, RI, 2005, pp.~1--15.

\bibitem{AngenentHulshofMatano}
{\sc S.~B. Angenent, J.~Hulshof, and H.~Matano}, {\em The radius of vanishing
  bubbles in equivariant harmonic map flow from {$D^2$} to {$S^2$}}, SIAM J.
  Math. Anal., 41 (2009), pp.~1121--1137.

\bibitem{BKP}
{\sc S.~Bartels, J.~Ko, and A.~Prohl}, {\em Numerical analysis of an explicit
  approximation scheme for the {L}andau-{L}ifshitz-{G}ilbert equation}, Math.
  Comp., 77 (2008), pp.~773--788.

\bibitem{BertschDalPassoHout}
{\sc M.~Bertsch, R.~Dal~Passo, and R.~van~der Hout}, {\em Nonuniqueness for the
  heat flow of harmonic maps on the disk}, Arch. Ration. Mech. Anal., 161
  (2002), pp.~93--112.

\bibitem{BertschHulshofHout}
{\sc M.~Bertsch, J.~Hulshof, and R.~van~der Hout}, {\em Energy concentration
  for 2-dimensional radially symmetric equivariant harmonic map heat flows},
  2010.
\newblock Preprint.

\bibitem{BertschPodioValente}
{\sc M.~Bertsch, P.~Podio~Guidugli, and V.~Valente}, {\em On the dynamics of
  deformable ferromagnets. {I}. {G}lobal weak solutions for soft ferromagnets
  at rest}, Ann. Mat. Pura Appl. (4), 179 (2001), pp.~331--360.

\bibitem{BBCH}
{\sc F.~Bethuel, H.~Brezis, B.~D. Coleman, and F.~H{\'e}lein}, {\em Bifurcation
  analysis of minimizing harmonic maps describing the equilibrium of nematic
  phases between cylinders}, Arch. Rational Mech. Anal., 118 (1992),
  pp.~149--168.

\bibitem{BW2}
{\sc C.~J. Budd and J.~F. Williams}, {\em Moving mesh generation using the
  parabolic monge-amp\`ere equation}, SIAM J. Sci. Comp., 31 (2009),
  pp.~3438--3465.

\bibitem{BuddWilliams}
\leavevmode\vrule height 2pt depth -1.6pt width 23pt, {\em How to adaptively
  resolve evolutionary singularities in differential equations with symmetry},
  J. Engrg. Math., 66 (2010), pp.~217--236.

\bibitem{CDY}
{\sc K.-C. Chang, W.~Y. Ding, and R.~Ye}, {\em Finite-time blow-up of the heat
  flow of harmonic maps from surfaces}, J. Differential Geom., 36 (1992),
  pp.~507--515.

\bibitem{CSU}
{\sc N.-H. Chang, J.~Shatah, and K.~Uhlenbeck}, {\em Schr\"odinger maps}, Comm.
  Pure Appl. Math., 53 (2000), pp.~590--602.

\bibitem{Coron}
{\sc J.-M. Coron}, {\em Nonuniqueness for the heat flow of harmonic maps}, Ann.
  Inst. H. Poincar\'e Anal. Non Lin\'eaire, 7 (1990), pp.~335--344.

\bibitem{Freire}
{\sc A.~Freire}, {\em Uniqueness for the harmonic map flow in two dimensions},
  Calc. Var. Partial Differential Equations, 3 (1995), pp.~95--105.

\bibitem{GuoHong}
{\sc B.~L. Guo and M.~C. Hong}, {\em The {L}andau-{L}ifshitz equation of the
  ferromagnetic spin chain and harmonic maps}, Calc. Var. Partial Differential
  Equations, 1 (1993), pp.~311--334.

\bibitem{GustKangTsai}
{\sc S.~Gustafson, K.~Kang, and T.-P. Tsai}, {\em Schr\"odinger flow near
  harmonic maps}, Comm. Pure Appl. Math., 60 (2007), pp.~463--499.

\bibitem{GKT2}
\leavevmode\vrule height 2pt depth -1.6pt width 23pt, {\em Asymptotic
  stability, concentration, and oscillation in harmonic map heat-flow,
  {L}andau-{L}ifshitz, and {S}chr\"odinger maps on {$\mathbb{R}^2$}}, Comm.
  Math. Phys., 300 (2010), pp.~205--242.

\bibitem{Hatcher}
{\sc A.~Hatcher}, {\em Algebraic topology}, Cambridge University Press,
  Cambridge, 2002.

\bibitem{HR}
{\sc W.~Huang, Y.~Ren, and R.~D. Russell}, {\em Moving mesh partial
  differential equations ({MMPDES}) based on the equidistribution principle},
  SIAM J. Numer. Anal., 31 (1994), pp.~709--730.

\bibitem{HS}
{\sc A.~Hubert and R.~Sch\"afer}, {\em Magnetic Domains: The Analysis of
  Magnetic Microstructures}, Springer, Berlin–Heidelberg–New York, 1998.

\bibitem{KP}
{\sc M.~Kru{\v{z}}{\'{\i}}k and A.~Prohl}, {\em Recent developments in the
  modeling, analysis, and numerics of ferromagnetism}, SIAM Rev., 48 (2006),
  pp.~439--483 (electronic).

\bibitem{Lemaire}
{\sc L.~Lemaire}, {\em Applications harmoniques de surfaces riemanniennes}, J.
  Differential Geom., 13 (1978), pp.~51--78.

\bibitem{MRR}
{\sc F.~Merle, P.~Rapha{\"e}l, and I.~Rodnianski}, {\em Blow up dynamics for
  smooth data equivariant solutions to the energy critical schr{\"o}dinger map
  problem}, 2011.
\newblock Preprint, arXiv:1102.4308.

\bibitem{RS}
{\sc P.~Rapha{\"e}l and R.~Schweyer}, {\em Stable blow up dynamics for the
  1-corotational energy critical harmonic heat flow}, 2011.
\newblock Preprint, arXiv:1106.0914.

\bibitem{VAGbook}
{\sc A.~A. Samarskii, V.~A. Galaktionov, S.~P. Kurdyumov, and A.~P. Mikhailov},
  {\em Blow-up in quasilinear parabolic equations}, vol.~19 of de Gruyter
  Expositions in Mathematics, Walter de Gruyter \& Co., Berlin, 1995.
\newblock Translated from the 1987 Russian original by Michael Grinfeld and
  revised by the authors.

\bibitem{Struwe2}
{\sc M.~Struwe}, {\em On the evolution of harmonic mappings of {R}iemannian
  surfaces}, Comment. Math. Helv., 60 (1985), pp.~558--581.

\bibitem{Struwe}
\leavevmode\vrule height 2pt depth -1.6pt width 23pt, {\em Geometric evolution
  problems}, in Nonlinear partial differential equations in differential
  geometry ({P}ark {C}ity, {UT}, 1992), vol.~2 of IAS/Park City Math. Ser.,
  Amer. Math. Soc., Providence, RI, 1996, pp.~257--339.

\bibitem{Topping}
{\sc P.~Topping}, {\em Reverse bubbling and nonuniqueness in the harmonic map
  flow}, Int. Math. Res. Not.,  (2002), pp.~505--520.

\bibitem{BHK}
{\sc J.~B. van~den Berg, J.~Hulshof, and J.~R. King}, {\em Formal asymptotics
  of bubbling in the harmonic map heat flow}, SIAM J. Appl. Math., 63 (2003),
  pp.~1682--1717 (electronic).

\bibitem{vanderHout}
{\sc R.~van~der Hout}, 2010.
\newblock personal communication.

\end{thebibliography}

\end{document}